# Fourier series (based) multiscale method for computational analysis in science and engineering:

# VII. Fourier series multiscale solution for elastic bending of beams on Pasternak foundations

Weiming Sun[a,*,+] and Zimao Zhang[b]

**Abstract:** Fourier series multiscale method, a concise and efficient analytical approach for multiscale computation, will be developed out of this series of papers. In the seventh paper, the usual structural analysis of beams on an elastic foundation is extended to a thorough multiscale analysis for a fourth order linear differential equation for transverse deflection of the beam, where general boundary conditions and a wide spectrum of model parameters are prescribed. For this purpose, the solution function is expressed as a linear combination of the boundary function and the internal function, to ensure the series expression obtained uniformly convergent and termwise differentiable up to fourth order. Meanwhile, the internal function corresponds to the particular solution, and the boundary function corresponds to the general solution which satisfies the homogeneous form of the governing differential equation. Since the general solution has appropriately interpreted the meaning of the differential equation, the spatial characteristics of the solution of the equation are expected to be better captured. With the boundary function and the internal function selected specifically as combination of the linearly independent homogeneous solutions of the differential equation, and one-dimensional half-range Fourier sine series over the solution interval, the Fourier series multiscale solution of the bending problem of a beam on the Pasternak foundation is derived. And then the convergence characteristics of the Fourier series multiscale solution are investigated with numerical examples, and the multiscale characteristics of the bending problem of a beam on the Pasternak foundation are demonstrated for a wide spectrum of model parameters.

---

[a]Department of Mathematics and Big Data, School of Artificial Intelligence, Jianghan University, Wuhan, 430056, China

[b]Department of Mechanics, School of Civil Engineering, Beijing Jiaotong University, Beijing, 100044, China

*Correspondence to: Weiming Sun, Department of Mathematics and Big Data, School of Artificial Intelligence, Jianghan University, Wuhan, 430056, China

[+] E-mail: xuxinenglish@hust.edu.cn

# 1. Introduction

The analysis of a beam on an elastic foundation is often performed in modern civil engineering. In this process, various beam theories, such as the classical Bernoulli-Euler beam theory [1-6] or the well-known Timoshenko beam theory [7] are adopted. Meanwhile, the foundation is usually idealized as a Winkler (one-parameter) model [1, 3, 6], a Pasternak (two-parameter) model [2, 4, 5], or even a Kerr (three-parameter) model [7]. Thanks to these simple mathematical formulations, they have been employed in a variety of problems and given satisfactory results in many practical situations [8]. However, note that all the above-mentioned studies are confined to structural analysis of building, geotechnical, airport runway, highway, and railroad structures, in which the model parameters of the foundations vary within a specified (narrow) range.

In the last two decades, multiscale problems in science and engineering have attracted considerable attentions [9]. The spatial, material and/or physical parameters in a multiscale problem, can vary easily by orders of magnitude, resulting in the appearance of boundary layers or other forms of local discontinuities with large gradients in the solution region. For example, the convection-diffusion-reaction equation takes a very simple form of the second order linear differential equation. However, when the model parameters change within a wide spectrum, it has varied physical behaviors with the presence of the exponential regime and the propagation regime [10-13]. And it becomes one of the most studied subjects in multiscale methods.

Compared to the convection-diffusion-reaction equation, the governing equation for elastic bending of Bernoulli-Euler beams on Pasternak foundations seems a little more complicated. Therefore, when the investigations are extended to a wide spectrum of model parameters, the specific structural analysis of engineering structures is expected to turn into a thorough multiscale analysis for a fourth order linear differential equation.

Nowadays, many multiscale methods have been developed for the multiscale problems, which include the stabilized finite element methods [14, 15], the bubble methods [16-19], the wavelet finite element methods [20-22], the meshless methods [23], the variational multiscale methods [24, 25], and so on. It is evident from their names that these multiscale methods are typically based on the traditional mesh-based or other discretization methods by making some suitable modifications, such as: use of stabilization terms, inclusion of different scale groups, decomposition of the solution into coarse and fine scale components, *etc*. Although these multiscale methods have been applied to a variety of boundary value problems, such as the convection-diffusion-reaction equation, we are still in the constant struggles with respect to the stability of the solution algorithms, proper selection of computational scales, robustness and effectiveness of the methods, balancing the span of scale groups and solution accuracy, high computational costs, low numerical accuracy for higher order derivatives of field variables, and so on.

It's worthy of notice that, a new type Fourier series method, namely the Fourier series method with supplementary terms, potentially represents a strategic shift from the existing framework towards better resolving multiscale problems. In this method, the solution function is expressed as a linear combination of a conventional Fourier series and some supplementary terms [26-44]. The supplementary terms, on the one hand, are purposely introduced to carry over the discontinuities potentially associated with the periodic extensions of the original solution function. As a result, the Fourier series actually corresponds to a periodic and sufficiently smooth residual function, and is hence uniformly convergent and termwise differentiable. On the other hand, the supplementary terms are only required to be sufficiently smooth over a compact interval (or domain), without regulating them to any particular forms. This implies that there is actually a large (theoretically, an infinite) number of possible choices for such basis functions in the process of implementing the Fourier series expansions. Specifically, if the supplementary terms are sought as the general solution of the differential equation, it can be expected that such general solution shall be able to appropriately interpret the meaning of the differential equation, and hence better capture the spatial characteristics of the solution in separate directions. In view of the multiscale capability of this solution method, we rename it the Fourier series (based) multiscale method.

Therefore, in the seventh paper of the series of researches on Fourier series multiscale method, we will reinvestigate the bending problem of a Bernoulli-Euler beam on the Pasternak foundation from a more systematic and general view of point. Firstly, we generalize this issue to the multiscale analysis of a fourth order linear differential equation for transverse deflection of the beam, where general boundary conditions and a wide spectrum of model parameters are prescribed. Secondly, the Fourier series multiscale method is developed for the analysis of the bending problem of a Bernoulli-Euler beam on the Pasternak foundation. As a routine task, the solution function, more exactly the transverse deflection, is decomposed into several constituents, such as the boundary function and the internal function. The internal function corresponds to the particular solution. And the boundary function corresponds to the general solution which satisfies the homogeneous form of the governing differential equation. And then, with all the constituents selected respectively as combination of the linearly independent homogeneous solutions of the differential equation, and one-dimensional half-range Fourier sine series over the solution interval, the Fourier series multiscale solution of the bending problem of a Bernoulli-Euler beam on the Pasternak foundation is derived. Accordingly, this paper begins with description of the problem. Detailed formulations related to the Fourier series multiscale method is then presented. Finally, convergence characteristics of the Fourier series multiscale solution are investigated with numerical examples, and the multiscale characteristics of the bending problem of a Bernoulli-Euler beam on the Pasternak foundation are demonstrated for a wide spectrum of model parameters.

## 2. Description of the problem

In this section, we make the following assumptions about the beam and the elastic foundation:

1. The beam is an elastic, straight beam with flexural rigidity $EI$ over the interval $[0, a]$, and its deformation and internal foreces, such as the beam slope, bending moment, and shear force, can be described by the beam deflection $w$.

2. For the beam, the upper surface is subjected to a transverse load $q$, and the lower surface only a reactive force $q_e$ of the foundation. A perfect contact exists between the beam and the foundation. In general, we have a biparametric representation [8] of the reactive force $q_e(x_1)$ of the foundation

$$q_e = kw - G_p \frac{d^2 w}{dx_1^2}, \tag{1}$$

where $k$ is the modulus of subgrade reaction of the foundation, $G_p$ is the shear modulus of the foundation.

According to Bernoulli-Euler beam theory, we can express the differential equation of equilibrium of transverse elastic bending of a beam resting on biparametric foundation as

$$\mathcal{L}_b w = q - q_e, \tag{2}$$

where the differential operator

$$\mathcal{L}_b = EI \frac{d}{dx_1^4}. \tag{3}$$

Meanwhile, we have the geometric equation

$$\theta = \frac{dw}{dx_1}, \tag{4}$$

and the physical equations

$$M = EI \frac{d^2w}{dx_1^2}, \tag{5}$$

$$Q = EI \frac{d^3w}{dx_1^3}, \tag{6}$$

where $\theta$, $M$ and $Q$ are the slope, bending moment and transverse shear force of the beam, respectively.

As an example, three kinds of boundary conditions are considered at the boundary (actually an end) $x_1 = a$ of the beam:

1. Generalized clamped boundary condition (C type boundary condition for short)

$$w(a) = \bar{w}_a, \quad \theta(a) = \bar{\theta}_a, \tag{7}$$

2. Generalized simply supported boundary condition (S type boundary condition for short)

$$w(a) = \bar{w}_a, \quad M(a) = \bar{M}_a, \tag{8}$$

3. Generalized free boundary condition (F type boundary condition for short)

$$M(a) = \bar{M}_a, \quad Q(a) = \bar{Q}_a, \tag{9}$$

where $\bar{w}_a$, $\bar{\theta}_a$, $\bar{M}_a$ and $\bar{Q}_a$ are the specified deflection, slope, bending moment and shear force at the end, respectively.

## 3. The Fourier series multiscale solution

The elastic bending of beams on biparametric foundations can be formulated as a fourth order ($2r_w = 4$) linear differential equation with constant coefficients of the transverse deflection of the beam. Moreover, this equation includes only the solution function and its even order derivatives. Therefore, the transverse deflection of the beam will be sought in the form of the one-dimensional half-range Fourier sine series multiscale solution as described below [45].

*3.1. The differential equation of beam deflection*

We substitute Eq. (1) in Eq. (2), and obtain the differential equation of the beam deflection

$$\mathcal{L}_{bf} w(x_1) = q(x_1), \tag{10}$$

where the differential operator is denoted by

$$\mathcal{L}_{bf} = EI \frac{d^4}{dx_1^4} - G_p \frac{d^2}{dx_1^2} + k. \tag{11}$$

## 3.2. Structural decomposition of the solution function

According to the Fourier series multiscale method for the fourth order linear differential equation with constant coefficients [45], we expand the solution function $w(x_1)$ in composite half-range Fourier sine series over the interval $[0, a]$. It will be decomposed into two parts:

$$w(x_1) = w_0(x_1) + w_1(x_1), \tag{12}$$

where $w_1(x_1)$ is the boundary function such that

$$w_1^{(2k_1)}(a) = w^{(2k_1)}(a), \quad w_1^{(2k_1)}(0) = w^{(2k_1)}(0), \quad k_1 = 0, 1, \tag{13}$$

and $w_0(x_1)$ is the internal function and naturally satisfies the sufficient conditions for 4 times term-by-term differentiation of the half-range Fourier sine series expansion of one-dimensional functions [46]

$$w_0^{(2k_1)}(a) = 0, \quad w_0^{(2k_1)}(0) = 0, \quad k_1 = 0, 1. \tag{14}$$

Substitute Eq. (12) into Eq. (10), and further, suppose that

$$\mathcal{L}_{bf} w_1 = 0, \tag{15}$$

and

$$\mathcal{L}_{bf} w_0 = q, \tag{16}$$

then we can accordingly decompose the solution of Eq. (10) into two parts, namely, the general solution and the particular solution. In this decomposition, the general solution corresponds to the boundary function $w_1(x_1)$ of the composite Fourier series, such that

$$\left. \begin{array}{l} \mathcal{L}_{bf} w_1 = 0 \\ w_1^{(2k_1)}(a) = w^{(2k_1)}(a), \quad w_1^{(2k_1)}(0) = w^{(2k_1)}(0), \quad k_1 = 0, 1 \end{array} \right\}, \tag{17}$$

and the particular solution corresponds to the internal function $w_0(x_1)$ of the composite Fourier series, such that

$$\left. \begin{array}{l} \mathcal{L}_{bf} w_0 = q \\ w_0^{(2k_1)}(a) = 0, \quad w_0^{(2k_1)}(0) = 0, \quad k_1 = 0, 1 \end{array} \right\}. \tag{18}$$

However, for the linear differential equation with constant coefficients, such as Eq. (10), the characteristic of the external load function $q$ has effects to some extent on the convergence and computational accuracy of the Fourier series multiscale solution. Under some specific conditions, we can decompose the external load function into two parts, namely, the coarse scale component $q_s$ and the fine scale component $q - q_s$, to which the responses of the differential equation are determined with different methods.

More specifically, we select a function $w_s(x_1)$ of a specific form (e.g. the algebraical polynomial) such that

$$\mathcal{L}_{bf} w_s = q_s. \tag{19}$$

Substitute it into Eq. (10) and then we have

$$\mathcal{L}_{bf}(w - w_s) = q - q_s. \tag{20}$$

Without loss of generality, suppose that the Fourier series multiscale solution of Eq. (20), is

$$w - w_s = w_0 + w_1, \tag{21}$$

accordingly, the solution of Eq. (10) becomes

$$w = w_0 + w_1 + w_s. \tag{22}$$

For brevity, we also call Eq. (22) the Fourier series multiscale solution of Eq. (10). Hence, the Fourier series multiscale solution includes not only the necessary terms, $w_0$ and $w_1$, but also the optional and supplementary term $w_s$. With the given definitions of the general solution and the particular solution of Eq. (10), we can further denote the function $w_s$ by the supplementary solution of Eq. (10).

It will be observed that reasonable choice of the supplementary solution provides an effective approach for improvement of convergence and computational accuracy of the Fourier series multiscale solution.

*3.3. The general solution*

Suppose that the homogeneous solution of Eq. (17.a) is given as
$$p_H(x_1) = \exp(\eta x_1), \tag{23}$$
where $\eta$ is an undetermined constant.

Substituting Eq. (23) into Eq. (17.a), we obtain the characteristic equation
$$EI\eta^4 - G_p\eta^2 + k = 0. \tag{24}$$

We define the comparative modulus of subgrade reaction of the foundation
$$G_{pr} = \frac{G_p a^2}{EI}, \tag{25}$$
and the comparative shear modulus of the foundation
$$k_r = \frac{ka^4}{EI}. \tag{26}$$

It is easy to verify that Eq. (24) has the following four distinct real roots when $\Delta_r = G_{pr}^2 - 4k_r > 0$,
$$\eta_1 = \alpha_1, \ \eta_2 = -\alpha_1, \ \eta_3 = \alpha_2, \ \eta_4 = -\alpha_2, \tag{27}$$
where $\alpha_1 = \frac{1}{a}\sqrt{\frac{1}{2}\left[G_{pr} + \sqrt{G_{pr}^2 - 4k_r}\right]}$, $\alpha_2 = \frac{1}{a}\sqrt{\frac{1}{2}\left[G_{pr} - \sqrt{G_{pr}^2 - 4k_r}\right]}$.

Eq. (24) has two distinct double real roots when $\Delta_r = G_{pr}^2 - 4k_r = 0$,
$$\eta_1 = \eta_2 = \alpha_3, \ \eta_3 = \eta_4 = -\alpha_3, \tag{28}$$
where $\alpha_3 = \frac{1}{a}\sqrt{\frac{1}{2}G_{pr}}$.

Eq. (24) has four distinct complex roots when $\Delta_r = G_{pr}^2 - 4k_r < 0$,
$$\eta_{1,2,3,4} = \pm\alpha_5 \pm i\alpha_6, \tag{29}$$
where $\alpha_5 > 0$, $\alpha_6 > 0$, $i = \sqrt{-1}$ and $\alpha_5 + i\alpha_6 = \frac{1}{a}\sqrt{\frac{1}{2}\left[G_{pr} + i\sqrt{4k_r - G_{pr}^2}\right]}$.

Further, according to the distribution of the characteristic roots of characteristic equation (24), we can determine the four linearly independent homogeneous solutions $p_{l,H}(x_1)$, $l = 1, 2, 3, 4$, as presented in table 1.

Table 1: Expressions for $p_{l,H}(x_1)$, $l = 1, 2, 3, 4$.

| | $G_{pr}^2 - 4k_r > 0$ | $G_{pr}^2 - 4k_r = 0$ | $G_{pr}^2 - 4k_r < 0$ |
|---|---|---|---|
| $p_{1,H}(x_1)$ | $\dfrac{\sinh(\alpha_1 x_1)}{\sinh(\alpha_1 a)}$ | $\dfrac{\sinh(\alpha_3 x_1)}{\sinh(\alpha_3 a)}$ | $\dfrac{\sinh(\alpha_5 x_1)\sin(\alpha_6 x_1)}{\sinh(\alpha_5 a)\sin(\alpha_6 a)}$ |
| $p_{2,H}(x_1)$ | $\dfrac{\sinh[\alpha_1(a - x_1)]}{\sinh(\alpha_1 a)}$ | $\dfrac{x_1 \sinh(\alpha_3 x_1)}{a \sinh(\alpha_3 a)}$ | $\dfrac{\sinh(\alpha_5 x_1)\sin[\alpha_6(a - x_1)]}{\sinh(\alpha_5 a)\sin(\alpha_6 a)}$ |
| $p_{3,H}(x_1)$ | $\dfrac{\sinh(\alpha_2 x_1)}{\sinh(\alpha_2 a)}$ | $\dfrac{\sinh[\alpha_3(a - x_1)]}{\sinh(\alpha_3 a)}$ | $\dfrac{\sinh[\alpha_5(a - x_1)]\sin(\alpha_6 x_1)}{\sinh(\alpha_5 a)\sin(\alpha_6 a)}$ |
| $p_{4,H}(x_1)$ | $\dfrac{\sinh[\alpha_2(a - x_1)]}{\sinh(\alpha_2 a)}$ | $\dfrac{(a - x_1)\sinh[\alpha_3(a - x_1)]}{a \sinh(\alpha_3 a)}$ | $\dfrac{\sinh[\alpha_5(a - x_1)]\sin[\alpha_6(a - x_1)]}{\sinh(\alpha_5 a)\sin(\alpha_6 a)}$ |

When the computational parameter $k_r$ is adjusted from 1.0 to $10^2$, $10^4$ and $10^6$ successively, the evolution of the homogeneous solutions $p_{l,H}(x_1)$, $l = 1, 2, 3, 4$, are displayed respectively in Figures 1-3.

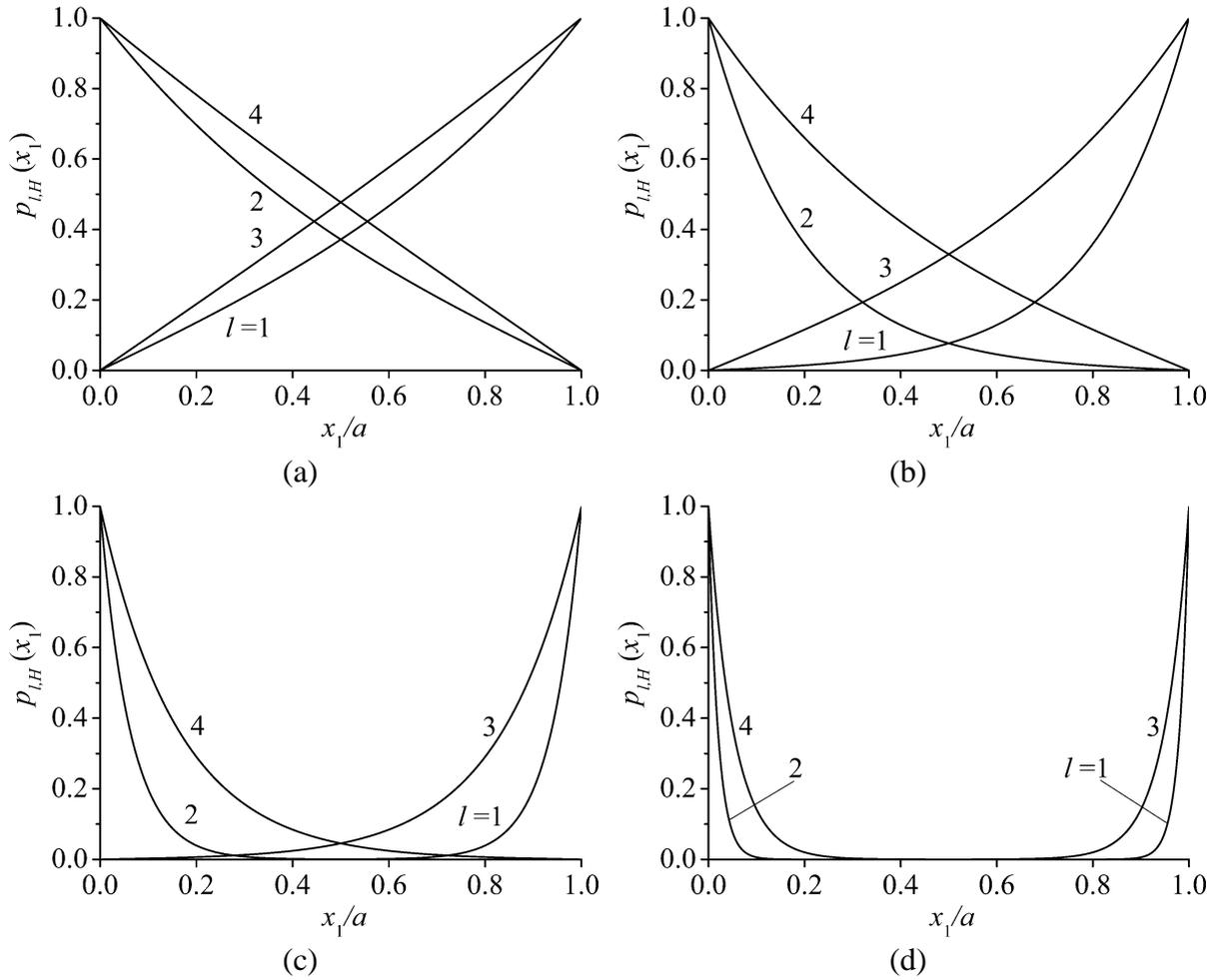

Figure 1: Homogeneous solutions with four different real roots:
(a) $k_r = 1.0$, $G_{pr} = 3\sqrt{k_r}$, (b) $k_r = 10^2$, $G_{pr} = 3\sqrt{k_r}$, (c) $k_r = 10^4$, $G_{pr} = 3\sqrt{k_r}$, (d) $k_r = 10^6$, $G_{pr} = 3\sqrt{k_r}$.

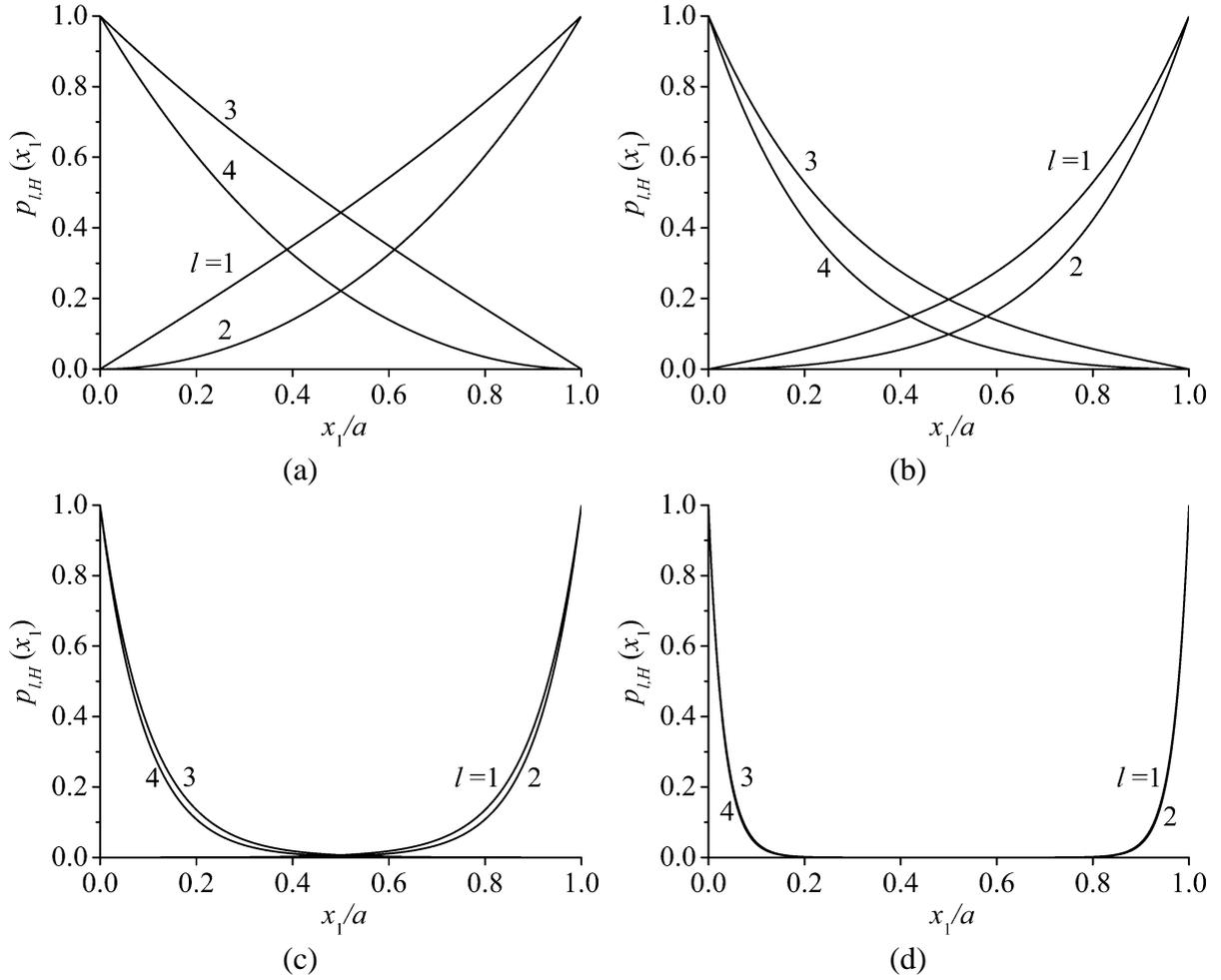

Figure 2: Homogeneous solutions with two real double roots:
(a) $k_r = 1.0$, $G_{pr} = 2\sqrt{k_r}$, (b) $k_r = 10^2$, $G_{pr} = 2\sqrt{k_r}$, (c) $k_r = 10^4$, $G_{pr} = 2\sqrt{k_r}$, (d) $k_r = 10^6$, $G_{pr} = 2\sqrt{k_r}$.

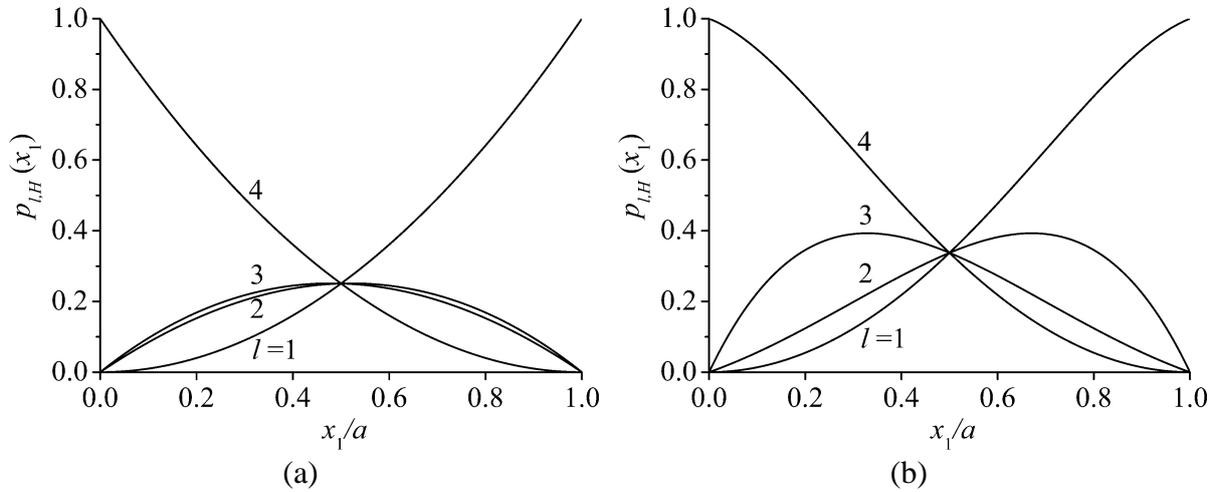

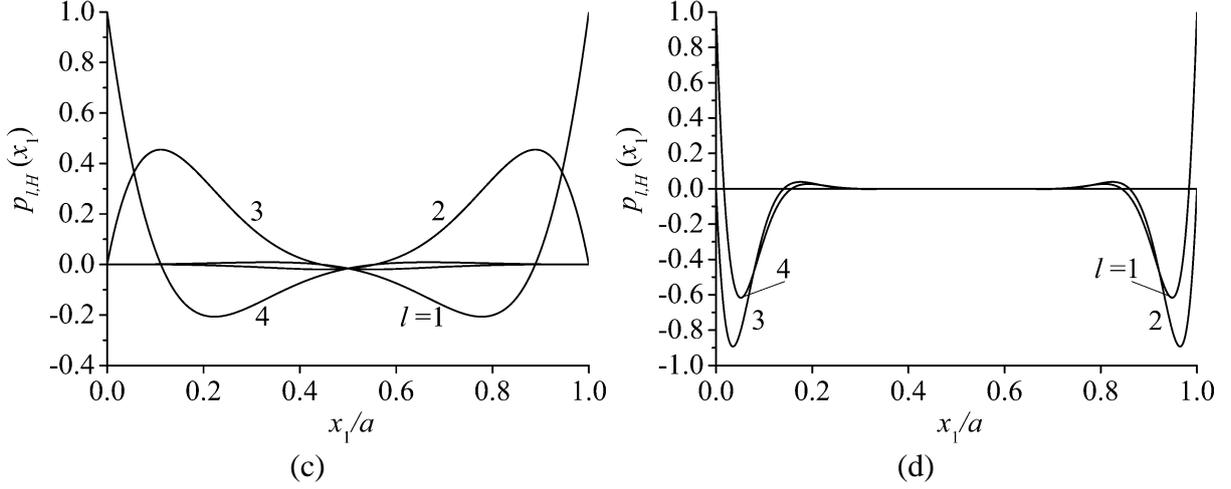

Figure 3: Homogeneous solutions with four different complex roots:
(a) $k_r = 1.0$, $G_{pr} = 0$, (b) $k_r = 10^2$, $G_{pr} = 0$, (c) $k_r = 10^4$, $G_{pr} = 0$, (d) $k_r = 10^6$, $G_{pr} = 0$.

We select the vector of functions
$$\mathbf{p}_1^T(x_1) = [p_{1,H}(x_1) \quad p_{2,H}(x_1) \quad p_{3,H}(x_1) \quad p_{4,H}(x_1)], \tag{30}$$
and for $k_1 = 0, 1, \cdots, 4$, we define the vector of derivatives of the selected functions
$$\mathbf{p}_1^{(k_1)T}(x_1) = [p_{1,H}^{(k_1)}(x_1) \quad p_{2,H}^{(k_1)}(x_1) \quad p_{3,H}^{(k_1)}(x_1) \quad p_{4,H}^{(k_1)}(x_1)]. \tag{31}$$

Further, let the vector of undetermined constants $\mathbf{a}_1^T = [G_{1,1} \quad G_{1,2} \quad G_{1,3} \quad G_{1,4}]$, then we can construct the boundary function as
$$w_1(x_1) = \mathbf{p}_1^T(x_1) \cdot \mathbf{a}_1. \tag{32}$$

Considering the necessary conditions as given in Eq. (17.b) for the boundary function $w_1(x_1)$, we have
$$\mathbf{R}_1 \mathbf{a}_1 = \mathbf{q}_1, \tag{33}$$
where
$$\mathbf{R}_1 = \begin{bmatrix} \mathbf{p}_1^{(0)T}(a) \\ \mathbf{p}_1^{(2)T}(a) \\ \mathbf{p}_1^{(0)T}(0) \\ \mathbf{p}_1^{(2)T}(0) \end{bmatrix}, \tag{34}$$
and the vector involving boundary values of the solution function $w(x_1)$ and its derivatives
$$\mathbf{q}_1^T = [w^{(0)}(a) \quad w^{(2)}(a) \quad w^{(0)}(0) \quad w^{(2)}(0)] = [w(a) \quad w^{(2)}(a) \quad w(0) \quad w^{(2)}(0)]. \tag{35}$$

Therefore, we have
$$\mathbf{a}_1 = \mathbf{R}_1^{-1} \mathbf{q}_1. \tag{36}$$
If we denote the vector of basis functions by
$$\mathbf{\Phi}_1^T(x_1) = \mathbf{p}_1^T(x_1) \cdot \mathbf{R}_1^{-1}, \tag{37}$$
then the general solution, corresponding to the boundary function $w_1(x_1)$, can finally be expressed as
$$w_1(x_1) = \mathbf{\Phi}_1^T(x_1) \cdot \mathbf{q}_1. \tag{38}$$

### 3.4. The supplementary solution

Let $N_1^s$ be a positive integer and $\Upsilon^s = \{x_{1,n_1}, n_1 = 1, 2, \cdots, N_1^s + 1\}$ be a set of interpolation points uniformly distributed on the interval $[0, a]$. Now we are to construct an interpolation algebraical polynomial $q_s(x_1)$ of the load function $q(x_1)$, such that the following interpolation conditions are satisfied

$$q_s(x_{1,n_1}) = q(x_{1,n_1}), \quad n_1 = 1, 2, \cdots, N_1^s + 1. \tag{39}$$

For this purpose, we select the vector of functions in the form of algebraical polynomials

$$\mathbf{p}_{qs}^T(x_1) = [p_{qs,1}(x_1) \quad p_{qs,2}(x_1) \quad \cdots \quad p_{qs, N_1^s+1}(x_1)], \tag{40}$$

where

$$p_{qs,j}(x_1) = (x_1/a)^{j-1}, \quad j = 1, 2, \cdots, N_1^s + 1, \tag{41}$$

and define the vector of undetermined constants

$$\mathbf{a}_{qs}^T = [H_{qs,1} \quad H_{qs,2} \quad \cdots \quad H_{qs, N_1^s+1}], \tag{42}$$

then the interpolation algebraical polynomial $q_s(x_1)$ can be expressed as

$$q_s(x_1) = \mathbf{p}_{qs}^T(x_1) \cdot \mathbf{a}_{qs}. \tag{43}$$

Combining the interpolation conditions in Eq. (39), we obtain

$$\mathbf{R}_{qs} \mathbf{a}_{qs} = \mathbf{q}_{qs}, \tag{44}$$

where

$$\mathbf{R}_{qs} = \begin{bmatrix} \mathbf{p}_{qs}^T(x_{1,1}) \\ \mathbf{p}_{qs}^T(x_{1,2}) \\ \vdots \\ \mathbf{p}_{qs}^T(x_{1,N_1^s+1}) \end{bmatrix}, \tag{45}$$

and the vector of values of the load function

$$\mathbf{q}_{qs}^T = [q(x_{1,1}) \quad q(x_{1,2}) \quad \cdots \quad q(x_{1,N_1^s+1})]. \tag{46}$$

Further, we consider the supplementary solution of Eq. (10) in the form of algebraical polynomials.

For $k \neq 0$, we select, without loss of generality, the vector of functions in the form of algebraical polynomials

$$\mathbf{p}_s^T(x_1) = [p_{s,1}(x_1) \quad p_{s,2}(x_1) \quad \cdots \quad p_{s, N_1^s+1}(x_1)], \tag{47}$$

where

$$p_{s,j}(x_1) = (x_1/a)^{j-1}, \quad j = 1, 2, \cdots, N_1^s + 1, \tag{48}$$

and define the vector of undetermined constants

$$\mathbf{a}_s^T = [G_{s,1} \quad G_{s,2} \quad \cdots \quad G_{s, N_1^s+1}], \tag{49}$$

thus, the supplementary solution of Eq. (10) can be expressed as

$$w_s(x_1) = \mathbf{p}_s^T(x_1) \cdot \mathbf{a}_s. \tag{50}$$

Substituting the supplementary solution into Eq. (19) we obtain

$$\mathbf{R}_s \mathbf{a}_s = \mathbf{a}_{qs}, \tag{51}$$

where the matrix

$$\mathbf{R}_s = [r_{s,ij}]_{i,j=1,2,\cdots,N_1^s+1}, \tag{52}$$

with

$$r_{s,ij} = \begin{cases} k, & \text{if } j=i, i=1,2,\cdots,N_1^s+1 \\ -G_p a^{-2} i(i+1), & \text{if } j=i+2, i=1,2,\cdots,N_1^s-1 \\ EI a^{-4} i(i+1)(i+2)(i+3), & \text{if } j=i+4, i=1,2,\cdots,N_1^s-3 \\ 0, & \text{otherwise} \end{cases}. \tag{53}$$

Hence

$$\mathbf{a}_s = \mathbf{R}_s^{-1} \mathbf{a}_{qs}. \tag{54}$$

Define the vector of basis functions

$$\boldsymbol{\Phi}_s^{\mathrm{T}}(x_1) = \mathbf{p}_s^{\mathrm{T}}(x_1) \cdot \mathbf{R}_s^{-1} \mathbf{R}_{qs}^{-1}. \tag{55}$$

Then the supplementary solution $w_s(x_1)$ can be expressed as

$$w_s(x_1) = \boldsymbol{\Phi}_s^{\mathrm{T}}(x_1) \cdot \mathbf{q}_{qs}. \tag{56}$$

In some occasions the supplementary solution $w_s(x_1)$ is not introduced into the Fourier series multiscale solution, $N_1^s = 0$ is denoted for convenience.

## 3.5. The particular solution

Let the error of the interpolation algebraical polynomial $q_s(x_1)$ relative to the load function $q(x_1)$ be

$$q_p(x_1) = q(x_1) - q_s(x_1), \tag{57}$$

then with the interpolation conditions we observe that

$$q_p(0) = 0, \quad q_p(a) = 0. \tag{58}$$

Expand it in a half-range Fourier sine series on the interval $[0, a]$, we have

$$q_p(x_1) = \sum_{m=1}^{\infty} V_{qp,2m} \sin(\alpha_m x_1), \tag{59}$$

where $\alpha_m = m\pi/a$, and $V_{qp,2m}$ are the Fourier coefficients of $q_p(x_1)$.

Denote

$$\varphi_{02m}(x_1) = \sin(\alpha_m x_1), \quad m = 1, 2, \cdots, \tag{60}$$

then we have:
for a nonnegative even integer $k_1$,

$$\varphi_{02m}^{(k_1)}(x_1) = (-1)^{k_1/2} \alpha_m^{k_1} \sin(\alpha_m x_1); \tag{61}$$

for a nonnegative odd integer $k_1$,

$$\varphi_{02m}^{(k_1)}(x_1) = (-1)^{(k_1-1)/2} \alpha_m^{k_1} \cos(\alpha_m x_1). \tag{62}$$

Further, $q_p(x_1)$ can be expressed in matrix form

$$q_p(x_1) = \boldsymbol{\Phi}_{02}^{\mathrm{T}}(x_1) \cdot \mathbf{q}_{qp,02}, \tag{63}$$

where the vector of trigonometric functions

$$\boldsymbol{\Phi}_{02}^{\mathrm{T}}(x_1) = [\varphi_{021}(x_1) \quad \varphi_{022}(x_1) \quad \cdots \quad \varphi_{02m}(x_1) \quad \cdots], \tag{64}$$

and the vector of Fourier coefficients

$$\mathbf{q}_{qp,02}^{\mathrm{T}} = [V_{qp,21} \quad V_{qp,22} \quad \cdots \quad V_{qp,2m} \quad \cdots]. \tag{65}$$

Meanwhile, suppose that Eq. (10) has the particular solution in the form of Fourier series
$$w_0(x_1) = \mathbf{\Phi}_{02}^T(x_1) \cdot \mathbf{q}_{02}, \tag{66}$$
where the vector of undetermined Fourier coefficients
$$\mathbf{q}_{02}^T = [V_{21} \quad V_{22} \quad \cdots \quad V_{2m} \quad \cdots]. \tag{67}$$
Substitute Eqs. (63) and (66) into Eq. (20), then we obtain
$$[\mathcal{L}_{bf}\mathbf{\Phi}_{02}^T(x_1)] \cdot \mathbf{q}_{02} = \mathbf{\Phi}_{02}^T(x_1) \cdot \mathbf{q}_{qp,02}, \tag{68}$$
which can be solved by the Fourier coefficient comparison method (FCCM).

For the particular solution $w_0(x_1)$, let $M$ be the number of truncated terms of the Fourier series, and compare the first $M$ Fourier coefficients on the both sides of Eq. (68) successively, then we have
$$(EI\alpha_m^4 + G_p\alpha_m^2 + k)V_{2m} = V_{qp,2m}, \quad m = 1, 2, \cdots, M. \tag{69}$$

By this means, the vector of undetermined Fourier coefficients $\mathbf{q}_{02}$ is obtained.

*3.6. Expression of the Fourier series multiscale solution*

Putting Eqs. (66), (38) and (56) together, we thus express the Fourier series multiscale solution for elastic bending of beams on biparametric foundations as
$$\begin{aligned} w(x_1) &= w_0(x_1) + w_1(x_1) + w_s(x_1) \\ &= \mathbf{\Phi}_{02}^T(x_1) \cdot \mathbf{q}_{02} + \mathbf{\Phi}_1^T(x_1) \cdot \mathbf{q}_1 + \mathbf{\Phi}_s^T(x_1) \cdot \mathbf{q}_{qs}, \end{aligned} \tag{70}$$
where the vectors of undetermined constants $\mathbf{q}_{02}$ and $\mathbf{q}_{qs}$ are related to the load function $q(x_1)$ and its corresponding interpolation algebraical polynomial $q_s(x_1)$, and are determined by Eqs. (69), (59) and Eq. (46) respectively. However, the vector of undetermined constants $\mathbf{q}_1$ is to be determined by taking the equation above to satisfy the prescribed boundary conditions.

*3.7. Equivalent transformation of the solution*

In the analysis of elastic bending of beams on biparametric foundations, it appears from Section 3.6 that the involved Fourier coefficients and expansion constants can be fully determined from forcing the solution to satisfy the boundary conditions. However, this leads to difficulties in the application of variational methods, where prior satisfaction of the displacement type boundary conditions is required. Therefore, in this section we perform the equivalent transformation for the change of primary undetermined constants.

For brevity, suppose that the displacement type boundary conditions are expressed by
$$\mathbf{q}_b = [w(0) \quad w^{(1)}(0) \quad w(a) \quad w^{(1)}(a)]^T = [w(0) \quad \theta(0) \quad w(a) \quad \theta(a)]^T. \tag{71}$$
Then substituting Eq. (70) into Eq. (71), we obtain
$$\mathbf{q}_1 = \mathbf{R}_f^{-1}\mathbf{q}_b - \mathbf{R}_f^{-1}\begin{bmatrix} \mathbf{\Phi}_{02}^T(0) \cdot \mathbf{q}_{02} + \mathbf{\Phi}_s^T(0) \cdot \mathbf{q}_{qs} \\ \mathbf{\Phi}_{02}^{(1)T}(0) \cdot \mathbf{q}_{02} + \mathbf{\Phi}_s^{(1)T}(0) \cdot \mathbf{q}_{qs} \\ \mathbf{\Phi}_{02}^T(a) \cdot \mathbf{q}_{02} + \mathbf{\Phi}_s^T(a) \cdot \mathbf{q}_{qs} \\ \mathbf{\Phi}_{02}^{(1)T}(a) \cdot \mathbf{q}_{02} + \mathbf{\Phi}_s^{(1)T}(a) \cdot \mathbf{q}_{qs} \end{bmatrix}, \tag{72}$$

where

$$\mathbf{R}_f = \begin{bmatrix} \mathbf{\Phi}_1^T(0) \\ \mathbf{\Phi}_1^{(1)T}(0) \\ \mathbf{\Phi}_1^T(a) \\ \mathbf{\Phi}_1^{(1)T}(a) \end{bmatrix}. \tag{73}$$

Plugging back into Eq. (70) gives

$$w(x_1) = \mathbf{\Phi}_1^T(x_1) \cdot \mathbf{R}_f^{-1} \mathbf{q}_b + \mathbf{\Phi}_{02}^T(x_1) \cdot \mathbf{q}_{02} + \mathbf{\Phi}_s^T(x_1) \cdot \mathbf{q}_{qs}$$

$$-\mathbf{\Phi}_1^T(x_1) \cdot \mathbf{R}_f^{-1} \begin{bmatrix} \mathbf{\Phi}_{02}^T(0) \cdot \mathbf{q}_{02} + \mathbf{\Phi}_s^T(0) \cdot \mathbf{q}_{qs} \\ \mathbf{\Phi}_{02}^{(1)T}(0) \cdot \mathbf{q}_{02} + \mathbf{\Phi}_s^{(1)T}(0) \cdot \mathbf{q}_{qs} \\ \mathbf{\Phi}_{02}^T(a) \cdot \mathbf{q}_{02} + \mathbf{\Phi}_s^T(a) \cdot \mathbf{q}_{qs} \\ \mathbf{\Phi}_{02}^{(1)T}(a) \cdot \mathbf{q}_{02} + \mathbf{\Phi}_s^{(1)T}(a) \cdot \mathbf{q}_{qs} \end{bmatrix}. \tag{74}$$

We denote the vector of functions

$$\mathbf{\Phi}_b^T(x_1) = \mathbf{\Phi}_1^T(x_1) \cdot \mathbf{R}_f^{-1}, \tag{75}$$

$$\mathbf{\Phi}_{0R}^T(x_1) = \mathbf{\Phi}_{02}^T(x_1) - \mathbf{\Phi}_b^T(x_1) \cdot \begin{bmatrix} \mathbf{\Phi}_{02}^T(0) \\ \mathbf{\Phi}_{02}^{(1)T}(0) \\ \mathbf{\Phi}_{02}^T(a) \\ \mathbf{\Phi}_{02}^{(1)T}(a) \end{bmatrix}, \tag{76}$$

$$\mathbf{\Phi}_{sR}^T(x_1) = \mathbf{\Phi}_s^T(x_1) - \mathbf{\Phi}_b^T(x_1) \cdot \begin{bmatrix} \mathbf{\Phi}_s^T(0) \\ \mathbf{\Phi}_s^{(1)T}(0) \\ \mathbf{\Phi}_s^T(a) \\ \mathbf{\Phi}_s^{(1)T}(a) \end{bmatrix}, \tag{77}$$

and the functions

$$w_{0R}(x_1) = \mathbf{\Phi}_{0R}^T(x_1) \cdot \mathbf{q}_{02}, \tag{78}$$

$$w_b(x_1) = \mathbf{\Phi}_b^T(x_1) \cdot \mathbf{q}_b, \tag{79}$$

$$w_{sR}(x_1) = \mathbf{\Phi}_{sR}^T(x_1) \cdot \mathbf{q}_{qs}, \tag{80}$$

then the Fourier series multiscale solution for elastic bending of beams on biparametric foundations can be rewritten as

$$w(x_1) = w_{0R}(x_1) + w_b(x_1) + w_{sR}(x_1). \tag{81}$$

Even though the solution given by Eq. (81) still looks similar to its origin Eq. (70), the displacement type boundary conditions have now been explicitly incorporated into the solution. In other words, the original coefficients associated with the boundary function, $w_1(x_1)$, have been transformed into a new set of coefficients which specifically depend upon the displacement type boundary conditions.

## 4. Variational Method

In addition to the usually used Fourier coefficient comparison method (FCCM), the variational method (VM), or more exactly, the minimum potential energy method, is adopted for the analysis of beams resting on foundations.

In this occassion, the supplementary solution $w_{sR}(x_1)$ is not introduced into the Fourier series multiscale solution. Therefore, in Eq. (81), we denote

$$\mathbf{q}_{qs} = \mathbf{0}, \tag{82}$$

and write the modified vector of basis functions $\mathbf{\Phi}_R^T$ and the modified vector of coefficients $\mathbf{q}_R^T$ respectively as

$$\mathbf{\Phi}_R^T = \begin{bmatrix} \mathbf{\Phi}_{0R}^T & \mathbf{\Phi}_b^T \end{bmatrix}, \tag{83}$$

$$\mathbf{q}_R^T = [\mathbf{q}_{02}^T \quad \mathbf{q}_b^T]. \tag{84}$$

And accordingly, Eq. (81) becomes

$$w = \mathbf{\Phi}_R^T \cdot \mathbf{q}_R. \tag{85}$$

With Eqs. (4)-(6), the transverse deflection, slope, bending moment and shear force of the beam can be expressed in matrix form

$$\begin{bmatrix} w \\ \theta \\ M \\ Q \end{bmatrix} = \mathbf{\Gamma}_R \cdot \mathbf{q}_R, \tag{86}$$

where the matrix

$$\mathbf{\Gamma}_R = \begin{bmatrix} \mathbf{\Gamma}_{R,w} \\ \mathbf{\Gamma}_{R,\theta} \\ \mathbf{\Gamma}_{R,M} \\ \mathbf{\Gamma}_{R,Q} \end{bmatrix} = \begin{bmatrix} \mathbf{\Phi}_R^T \\ \mathbf{\Phi}_R^{(1)T} \\ EI\mathbf{\Phi}_R^{(2)T} \\ EI\mathbf{\Phi}_R^{(3)T} \end{bmatrix}. \tag{87}$$

For the elastic system consisting of beam and biparametric foundation, we can formulate the energy of the system as

$$\Pi = \Pi_b + \Pi_f + \Pi_q + \Pi_\sigma, \tag{88}$$

where $\Pi_b$ is the elastic potential energy of beam, $\Pi_f$ is the elastic potential energy of biparametric foundation, $\Pi_q$ is the potential energy of transverse load distributed over the beam, and $\Pi_\sigma$ is the total potential energy of bending moment, and transverse shear force applied on the stress boundaries of the beam.

Specifically, the elastic potential energy of the beam is

$$\Pi_b = \frac{1}{2EI}\int_0^a M^2(x_1)dx_1. \tag{89}$$

The elastic potential energy of biparametric foundation is

$$\Pi_f = \frac{1}{2}\int_0^a [kw^2(x_1) + G_p(\frac{dw}{dx_1})^2]dx_1. \tag{90}$$

The potential energy of transverse load distributed over the beam is

$$\Pi_q = -\int_0^a qw(x_1)dx_1. \tag{91}$$

If all ends of the beam are prescribed as stress boundaries and there are bending moments and transverse shear forces applied at the ends, the corresponding potential energy is

$$\Pi_\sigma = -\overline{Q}_a w(a) - \overline{M}_a \theta(a) + \overline{Q}_0 w(0) + \overline{M}_0 \theta(0). \tag{92}$$

where $\overline{M}_0$ and $\overline{Q}_0$ are the specified bending moment and shear force at the end $x_1 = 0$.

Substituting Eqs. (89)-(92) in Eq. (88), we obtain

$$\Pi = \frac{1}{2}\mathbf{q}_R^T \mathbf{K}_{bf} \mathbf{q}_R - \mathbf{q}_R^T \mathbf{Q}_{bf}, \tag{93}$$

where the stiffness matrix

$$\mathbf{K}_{bf} = \int_0^a [\frac{1}{EI}\mathbf{\Gamma}_{R,M}^T \mathbf{\Gamma}_{R,M} + k\mathbf{\Gamma}_{R,w}^T \mathbf{\Gamma}_{R,w} + G_p \mathbf{\Gamma}_{R,\theta}^T \mathbf{\Gamma}_{R,\theta}]dx_1 \tag{94}$$

and the vector of the resultant forces

$$\mathbf{Q}_{bf} = \int_0^a q\mathbf{\Gamma}_{R,w}^T dx_1 + \overline{Q}_a \mathbf{\Gamma}_{R,w}^T(a) + \overline{M}_a \mathbf{\Gamma}_{R,\theta}^T(a) - \overline{Q}_0 \mathbf{\Gamma}_{R,w}^T(0) - \overline{M}_0 \mathbf{\Gamma}_{R,\theta}^T(0). \tag{95}$$

The corresponding stationary condition finally becomes

$$\mathbf{K}_{bf} \mathbf{q}_R = \mathbf{Q}_{bf}. \tag{96}$$

## 5. Numerical examples

As shown in Table 2, we have implemented the Fourier series multiscale method for elastic bending of beams on biparametric foundations with wide range of computational parameters and boundary conditions, and presented an analysis of the effects both of the order of interpolation algebraical polynomial of the load function and the derivation methods for discrete equations on the convergence characteristics and approximation accuracy of the Fourier series multiscale solution. By this means, we optimize the settings of computational schemes of the Fourier series multiscale solution of the elastic bending of beams on biparametric foundations. And with the optimized computational schemes, we investigate the multiscale characteristics of beams resting on biparametric foundations for varied computational parameters.

Table 2: Computational schemes of elastic bending of beams on biparametric foundations.

| Configuration of problem | | Computational scheme | |
|---|---|---|---|
| Boundary condition | Computational parameter $(k_r, G_{pr})$ | Order of interpolation algebraical polynomial of load function | Derivation method for discrete equations |
| CC | $(10^6, 0)$ | $N_1^s = 0$ | FCCM |
| SS | $(10^6, 1000)$ | $N_1^s = 1$ | VM |
| CF | $(10^6, 2000)$ | $N_1^s = 2$ | |
| | $(10^6, 3000)$ | | |

*5.1. Convergence characteristics*

As shown in Table 3, we perform three set of convergence comparison experiments in series for detailed analysis of the influences of some possible factors on convergence characteristics of the Fourier series multiscale solution. These possible factors are extended to involve in the governing differential equation for elastic bending of beams on biparametric foundations (computational parameters and boundary conditions) and its Fourier series multiscale solution (the order of the interpolation algebraical polynomial of the load function and the derivation method for discrete equations). We adopt the computational scheme specified in the first column of Table 3 as reference, and then adjustment of the single factors (see the corresponding rows in Table 3) leads to the three set of convergence comparison experiments. For instance, in the first set of experiments, we start out with the reference computational scheme, and change the order of the interpolation algebraical polynomial of the load function from zero to one and two successively. In the second set of experiments, we start out with the reference computational scheme, and introduce the variational method (VM) as the second derivation method for discrete equations, and then adjust the computational parameter $G_{pr}$ of the foundation from 0 to 1000, 2000 and 3000 successively. Finally, based on the reference computational scheme, we perform the third set of experiments by introducing the variational method as the second derivation method for discrete equations and adjusting boundary conditions from generalized clamped boundary condition (CC) to generalized simply supported boundary condition (SS) and generalized clamped-free boundary condition (CF) successively.

In the performed convergence comparison experiments, the load function remains unchanged, with the following form of the third order algebraical polynomial

$$q(x_1) = EIa^{-3}[10^3 + 2\times10^3 \frac{x_1}{a} + 5\times10^3 (\frac{x_1}{a})^2 + 10^4 (\frac{x_1}{a})^3], \quad x_1 \in [0,a], \tag{97}$$

Therefore, when we adjust the order of the interpolation algebraical polynomial of the load function to three, it follows that the vector of Fourier coefficients

$$\mathbf{q}_{02} = \mathbf{0}, \tag{98}$$

and accordingly the Fourier series multiscale solution given by Eq. (70), or equivalently Eq. (81), yields the exact solution of the elastic bending of beams on biparametric foundations.

Table 3: Convergence comparison experiments for elastic bending of beams on biparametric foundations.

| Numerical experiment | No. | Boundary condition | Computational parameter ($k_r$, $G_{pr}$) | Order of interpolation algebraical polynomial of load function | Derivation method for discrete equations |
|---|---|---|---|---|---|
| 1 | a | CC | ($10^6$, 0) | $N_1^s = 0$ | FCCM |
|   | b |    |              | $N_1^s = 1$ |      |
|   | c |    |              | $N_1^s = 2$ |      |
| 2 | a | CC | ($10^6$, 0)    | $N_1^s = 0$ | FCCM |
|   | b |    | ($10^6$, 1000) |             |      |
|   | c |    | ($10^6$, 2000) |             | VM   |
|   | d |    | ($10^6$, 3000) |             |      |
| 3 | a | CC | ($10^6$, 0)    | $N_1^s = 0$ | FCCM |
|   | b | SS |                |             | VM   |
|   | c | CF |                |             |      |

1. General convergence characteristics

As to the reference computational scheme given in Table 3, we truncate the composite Fourier series of the Fourier series multiscale solution successively with the first 2, 3, 5, 10, 20, 30 and 40 terms and then we compute the internal approximation errors and boundary approximation errors (see [47]) of the deflection $w(x_1)$, slope $\theta(x_1)$, and bending moment $M(x_1)$ of the beam. Some of the results are presented in Figure 4.

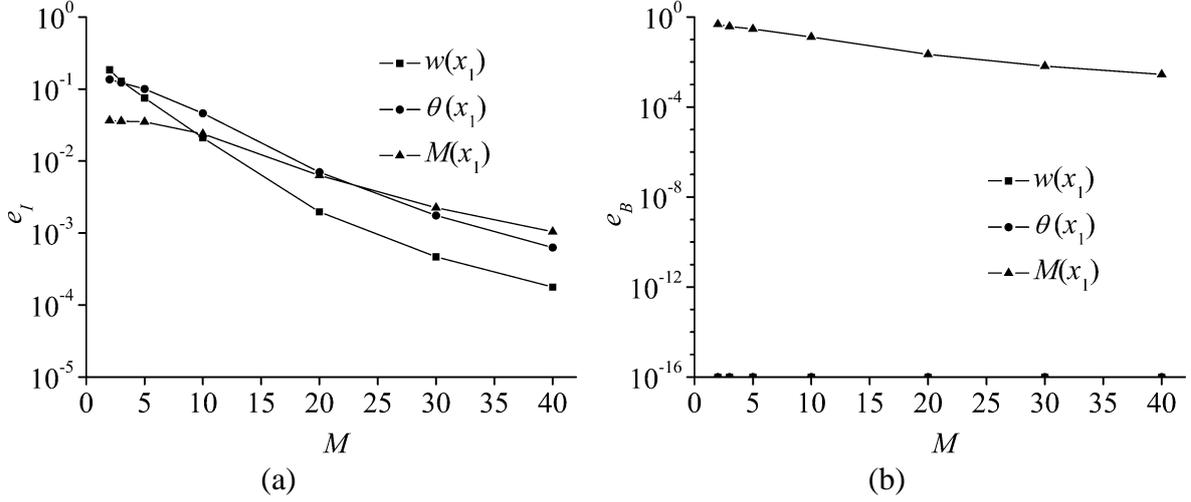

Figure 4: Convergence characteristics of the Fourier series multiscale solution for the elastic bending of a beam resting on the biparametric foundation:
(a) $e_I$-$M$ curves, (b) $e_B$-$M$ curves.

With the increase of the number of truncated terms of the composite Fourier series, the evolution of computed values of the indexes of approximation errors is exhibited as the following:

a. The internal approximation errors of $w(x_1)$ decrease rapidly as the number of truncated terms of the composite Fourier series increases. When the number of truncated terms equals to 10, the approximation error is already about 2.0E-2 and shows a trend of sustained and rapid decrease. And when the number of truncated terms equals to 40, the approximation error decreases further by 2 orders of magnitude.

b. The boundary approximation errors of $w(x_1)$ remain extremely small as the number of truncated terms of the composite Fourier series increases. The approximation errors are all below 1.0E-16.

c. The internal approximation errors of $\theta(x_1)$ decrease rapidly as the number of truncated terms of the composite Fourier series increases. When the number of truncated terms equals to 10, the approximation error is already about 5.0E-2 and shows a trend of sustained decrease. And when the number of truncated terms equals to 40, the approximation error decreases further by 2 orders of magnitude.

d. The boundary approximation errors of $\theta(x_1)$ remain extremely small as the number of truncated terms of the composite Fourier series increases. The approximation errors are all below 1.0E-16.

e. The internal approximation errors of $M(x_1)$ decrease rapidly as the number of truncated terms of the composite Fourier series increases. When the number of truncated terms equals to 10, the approximation error is already about 2.0E-2 and shows a trend of sustained decrease. And when the number of truncated terms equals to 40, the approximation error decreases further by 1 order of magnitude.

f. The boundary approximation errors of $M(x_1)$ decrease slowly as the number of truncated terms of the composite Fourier series increases. When the number of truncated terms equals to 10, the approximation error is already about 1.0E-1 and shows a trend of sustained decrease. And when the number of truncated terms equals to 40, the approximation error decreases further by 2 orders of magnitude.

Therefore, the convergence characteristics of the Fourier series multiscale solution for elastic bending of beams on biparametric foundations are summarized: the composite Fourier series of the deflection $w(x_1)$, slope $\theta(x_1)$, and bending moment $M(x_1)$ of the beam converge well not only within the solution interval, but also on the boundary.

2. Convergence comparison experiment 1

When the order of the interpolation algebraical polynomial of the load function changes from zero to one and two successively, the evolution of convergence characteristics of the Fourier series multiscale solution is illustrated in Figure 5. We make a brief analysis as the following:

a. The convergence of the composite Fourier series of $w(x_1)$, $\theta(x_1)$ and $M(x_1)$ remain almost unchanged within the solution interval, namely the composite Fourier series all converge well within the solution interval. However, when the number of truncated terms of the composite Fourier series equals to 40, the corresponding approximation errors decrease respectively by 4, 3 to 4, and 4 orders of magnitude.

b. The convergence of the composite Fourier series of $w(x_1)$ and $\theta(x_1)$ remain almost unchanged on the boundary of the solution interval, namely the composite Fourier series all converge well on the boundary. Moreover, the approximation errors of $w(x_1)$ and $\theta(x_1)$ remain extremely small (all below 1.0E-16) as the number of truncated terms of the composite Fourier series increases.

c. The convergence of the composite Fourier series of $M(x_1)$ remain almost unchanged on the boundary of the solution interval, namely the composite Fourier series all converge well on the boundary. When the number of truncated terms of the composite Fourier series equals to 40, the corresponding approximation errors all decrease by 4 orders of magnitude.

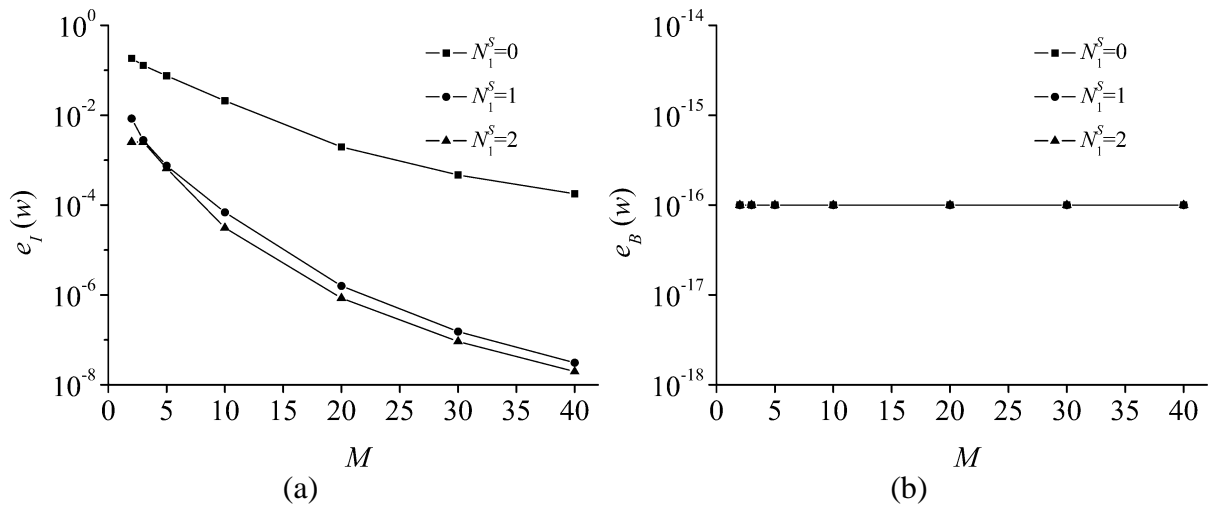

(a)          (b)

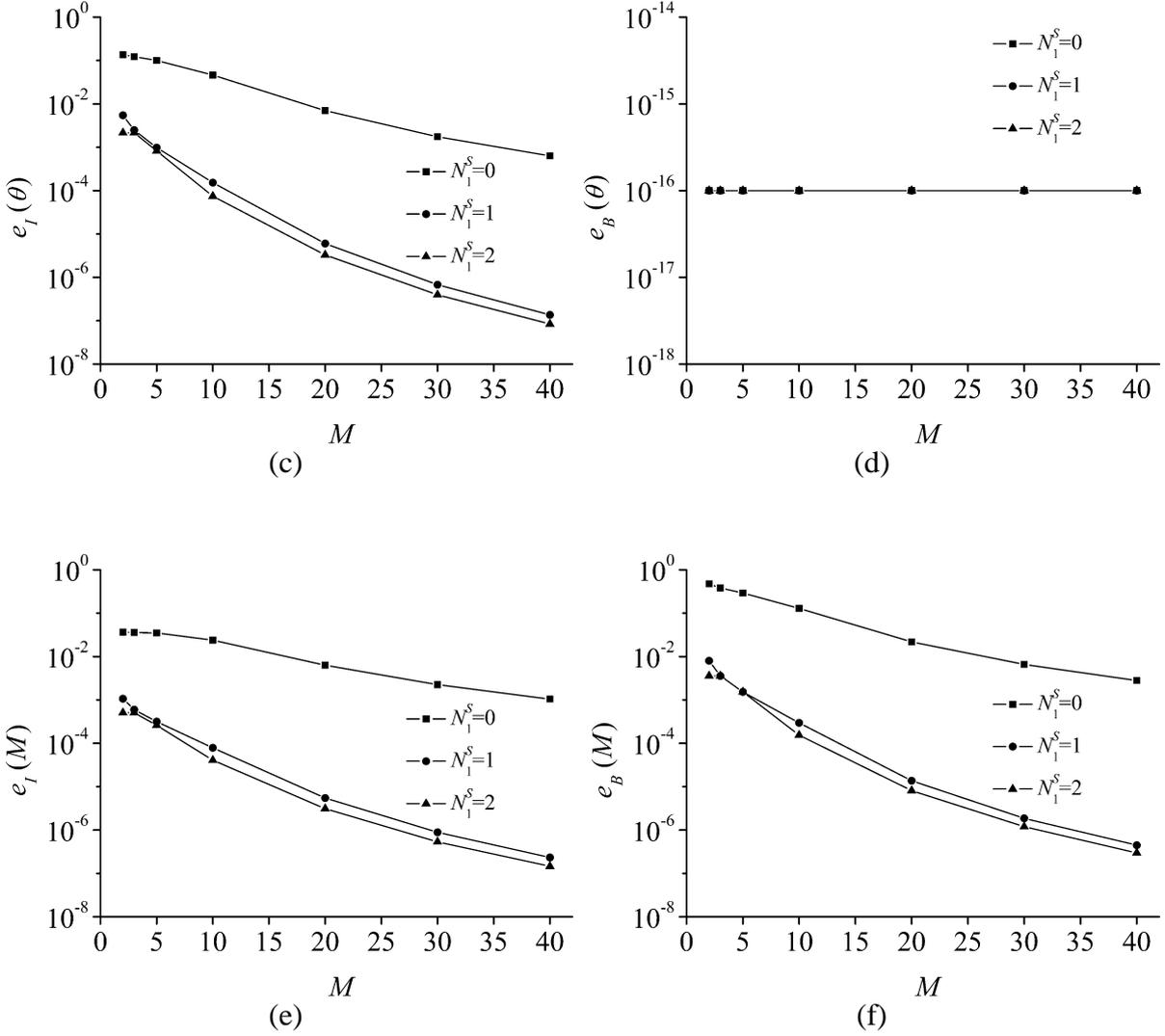

Figure 5: Convergence comparison of the Fourier series multiscale solutions with different orders of transverse load interpolation polynomials:
(a) $e_I(w)$-$M$ curves, (b) $e_B(w)$-$M$ curves, (c) $e_I(\theta)$-$M$ curves, (d) $e_B(\theta)$-$M$ curves, (e) $e_I(M)$-$M$ curves, (f) $e_B(M)$-$M$ curves.

3. Convergence comparison experiment 2

When one computational parameter of the foundation, $k_r$, is specified as $10^6$, and the other computational parameter $G_{pr}$ changes from 0 to 1000, 2000 and 3000 successively, the evolution of convergence characteristics of the Fourier series multiscale solution obtained by the Fourier coefficient comparison method (FCCM) is illustrated in Figure 6. We make a brief analysis as the following:

a. The convergence of the composite Fourier series of $w(x_1)$, $\theta(x_1)$ and $M(x_1)$ remain almost unchanged within the solution interval, namely the composite Fourier series all converge well within the solution interval. Meanwhile, when the number of truncated terms of the composite Fourier series equals to 40, the corresponding approximation errors are of no difference in quantity.

b. The convergence of the composite Fourier series of $w(x_1)$ and $\theta(x_1)$ remain almost

unchanged on the boundary of the solution interval, namely the composite Fourier series all converge well on the boundary. Moreover, the approximation errors of $w(x_1)$ and $\theta(x_1)$ remain extremely small (all below 1.0E-16) as the number of truncated terms of the composite Fourier series increases.

c. The convergence of the composite Fourier series of $M(x_1)$ remain almost unchanged on the boundary of the solution interval, namely the composite Fourier series all converge well on the boundary. And, when the number of truncated terms of the composite Fourier series equals to 40, the corresponding approximation errors of $M(x_1)$ are of no marked difference in quantity.

As a comparison, the evolution of convergence characteristics of the Fourier series multiscale solution obtained by the variational method (VM) is illustrated in Figure 7. And we come to the same conclusion.

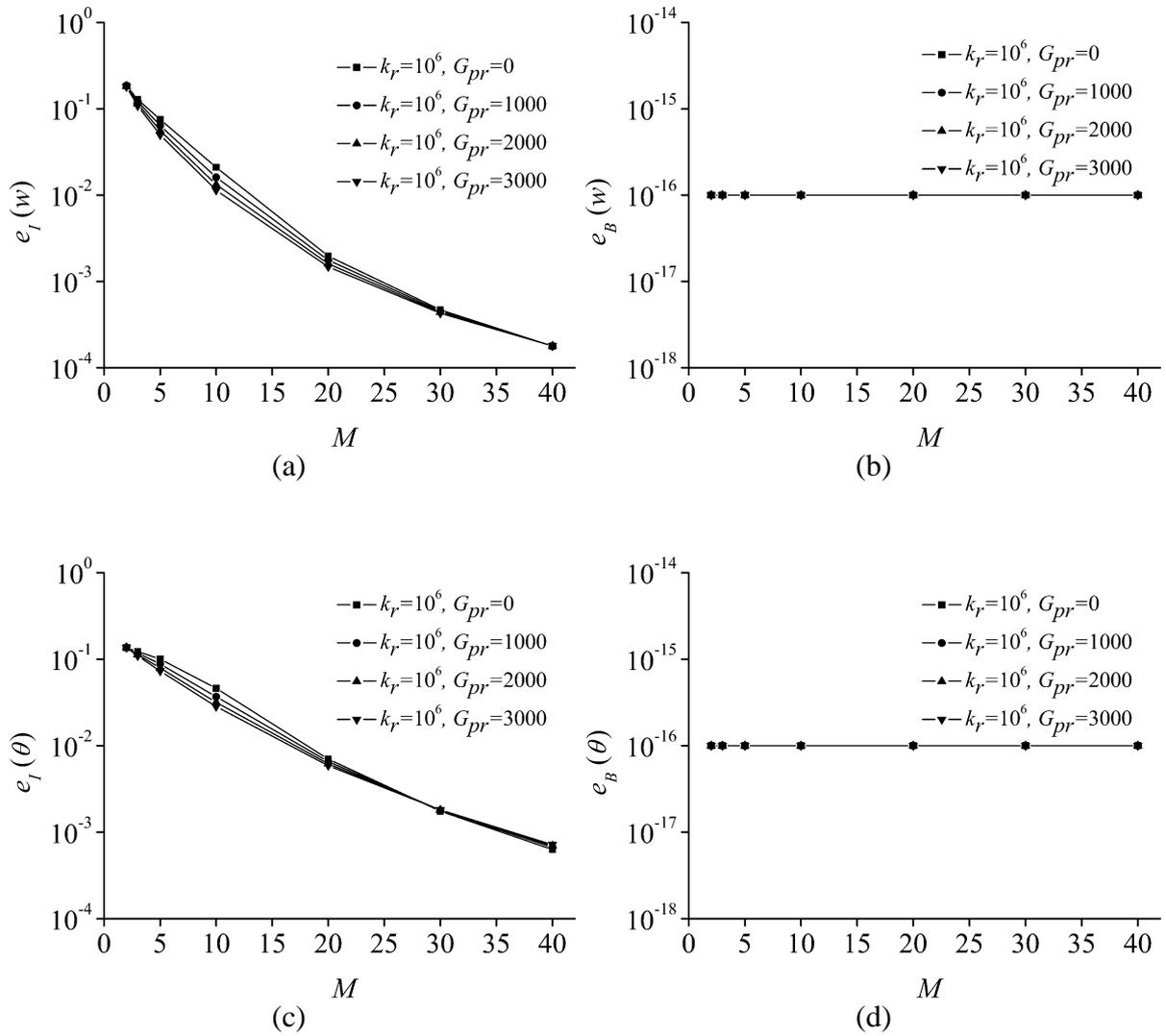

(a)

(b)

(c)

(d)

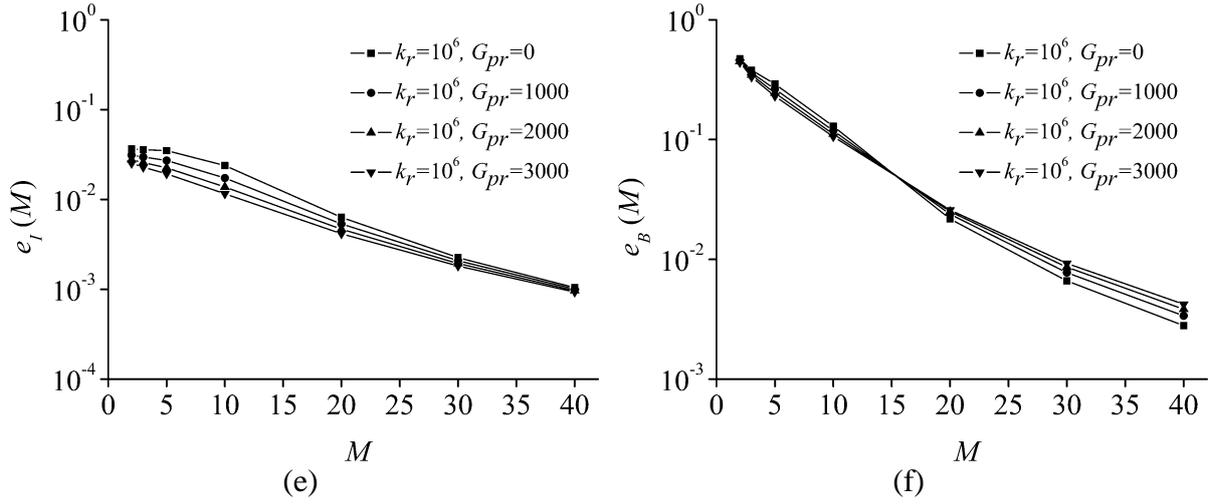

(e)   (f)

Figure 6: Convergence comparison of the Fourier series multiscale solutions obtained by FCCM, and with different computational parameters $k_r$ and $G_{pr}$:

(a) $e_I(w)$-$M$ curves, (b) $e_B(w)$-$M$ curves, (c) $e_I(\theta)$-$M$ curves, (d) $e_B(\theta)$-$M$ curves, (e) $e_I(M)$-$M$ curves, (f) $e_B(M)$-$M$ curves.

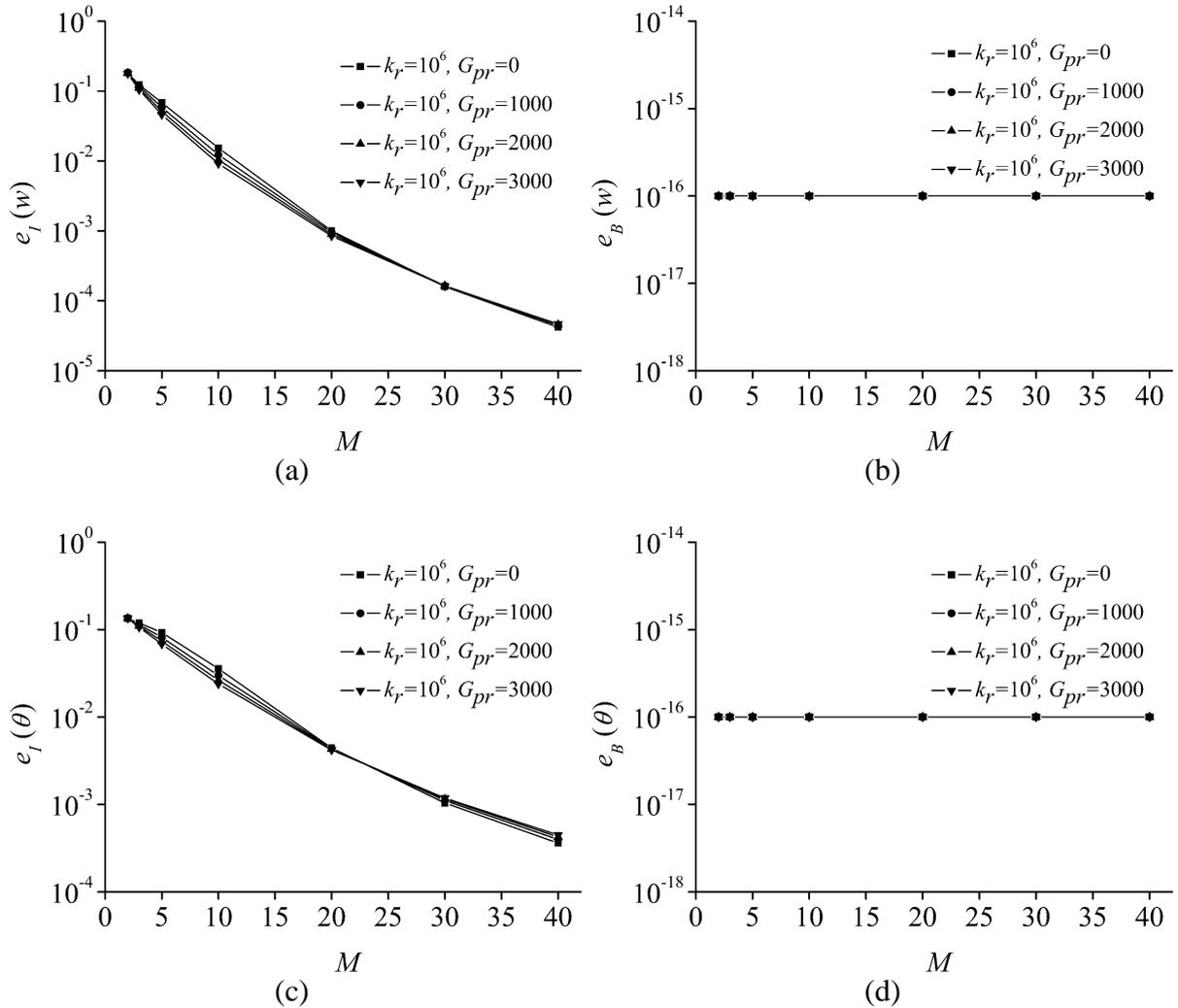

(a)   (b)

(c)   (d)

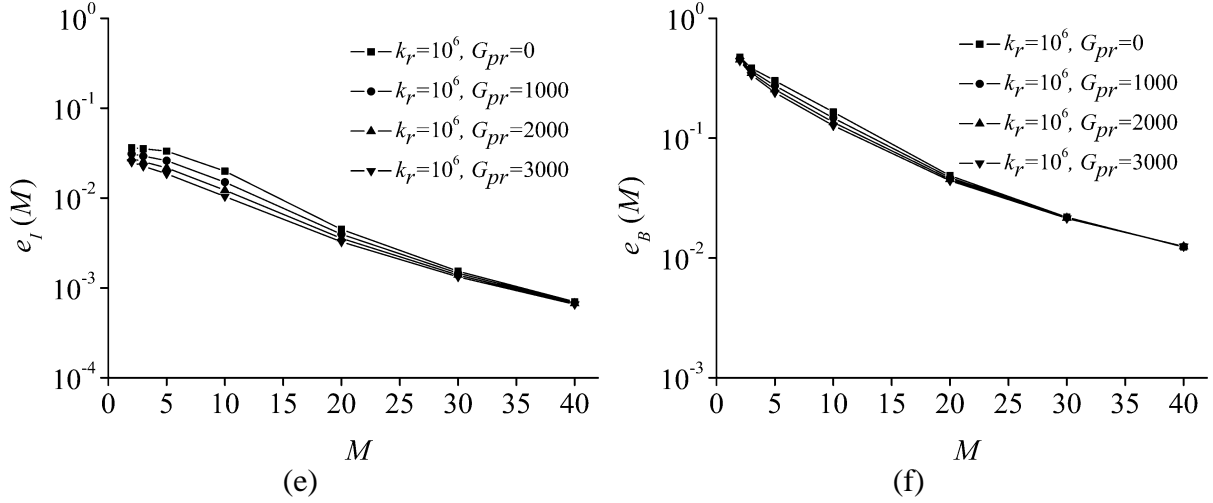

(e)                                         (f)

Figure 7: Convergence comparison of the Fourier series multiscale solutions obtained by VM, and with different computational parameters $k_r$ and $G_{pr}$:

(a) $e_I(w)$-$M$ curves, (b) $e_B(w)$-$M$ curves, (c) $e_I(\theta)$-$M$ curves, (d) $e_B(\theta)$-$M$ curves, (e) $e_I(M)$-$M$ curves, (f) $e_B(M)$-$M$ curves.

4. Convergence comparison experiment 3

When the boundary condition of elastic bending of beams on biparametric foundations changes from CC boundary condition to SS boundary condition and CF boundary condition, the evolution of convergence characteristics of the Fourier series multiscale solutions obtained by the Fourier coefficient comparison method (FCCM) is illustrated in Figure 8. We make a brief analysis as the following:

a. The convergence of the composite Fourier series of $w(x_1)$ remain almost unchanged within the solution interval, namely the composite Fourier series converge well within the solution interval. However, when the number of truncated terms of the composite Fourier series equals to 40, the corresponding approximation errors respectively increase by 0.5 orders of magnitude, and decrease by 1.5 orders of magnitude.

b. The convergence of the composite Fourier series of $\theta(x_1)$ and $M(x_1)$ within the solution interval varies with the boundary condition. For the CC boundary condition and SS boundary condition, the composite Fourier series converge well. And when the number of truncated terms of the composite Fourier series equals to 40, the corresponding approximation errors increase slightly. However, for the CF boundary condition, the convergence of the composite Fourier series of $\theta(x_1)$ and $M(x_1)$ is worsened greatly within the solution interval. when the number of truncated terms of the composite Fourier series equals to 40, the corresponding approximation errors increase respectively by 3 and 2 orders of magnitude.

c. The convergence of the composite Fourier series of $w(x_1)$, $\theta(x_1)$ and $M(x_1)$ on the boundary of the solution interval varies considerably with the boundary condition. For the CC boundary condition and SS boundary condition, the composite Fourier series converge well. More exactly, under these two boundary conditions, the approximation errors of $w(x_1)$, $\theta(x_1)$ and $M(x_1)$ remain extremely small (all below 1.0E-16) as the number of truncated terms of the composite Fourier series increases; or when the number of truncated terms of the composite Fourier series equals to 40, the corresponding approximation errors decrease to 3.0E-3. As to the CF boundary condition, the convergence of the composite Fourier series of $w(x_1)$, $\theta(x_1)$ is worsened considerably on the boundary of the solution interval. When the number of truncated terms of the composite Fourier series equals to 40, the corresponding approximation errors are respectively up to 0.11 and 1.13.

As a comparison, the evolution of convergence characteristics of the Fourier series multiscale solution obtained by the variational method (VM) is illustrated in Figure 9. We also make a brief analysis as the following:

a. The convergence of the composite Fourier series of $w(x_1)$, $\theta(x_1)$ and $M(x_1)$ remain almost unchanged within the solution interval, namely the composite Fourier series all converge well within the solution interval. Meanwhile, when the number of truncated terms of the composite Fourier series equals to 40, the corresponding approximation errors are of no difference in quantity, or increase respectively by 0-1 and 1-2 orders of magnitude.

b. The convergence of the composite Fourier series of $w(x_1)$, $\theta(x_1)$ and $M(x_1)$ remain almost unchanged on the boundary of the solution interval, namely the composite Fourier series all converge well on the boundary. Moreover, the approximation errors of $w(x_1)$ and $\theta(x_1)$ remain extremely small (all below 2.0E-7) as the number of truncated terms of the composite Fourier series increases. Moreover, when the number of truncated terms of the composite Fourier series equals to 40, the corresponding approximation errors of $M(x_1)$ are respectively 1.2E-2, 3.5E-2 and 1.8E-1 for the three boundary conditions.

Finally, when the number of truncated terms of the composite Fourier series increases, the corresponding curves for $w(x_1)$, $\theta(x_1)$ and $M(x_1)$ are presented respectively for the CC boundary condition in Figure 10 and for the CF boundary condition in Figure 11. And the advantage of the vaiational method (VM) has been fully demonstrated.

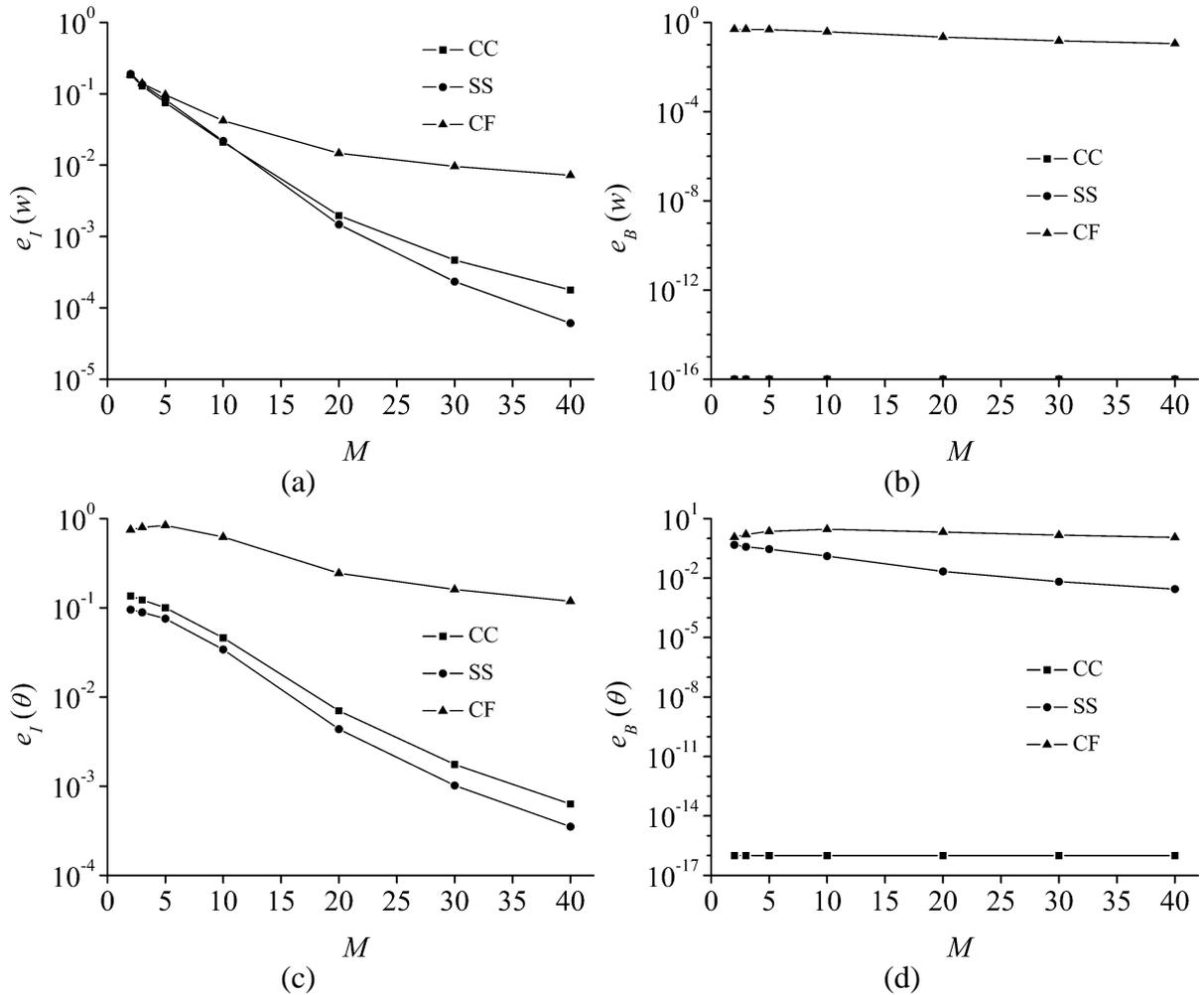

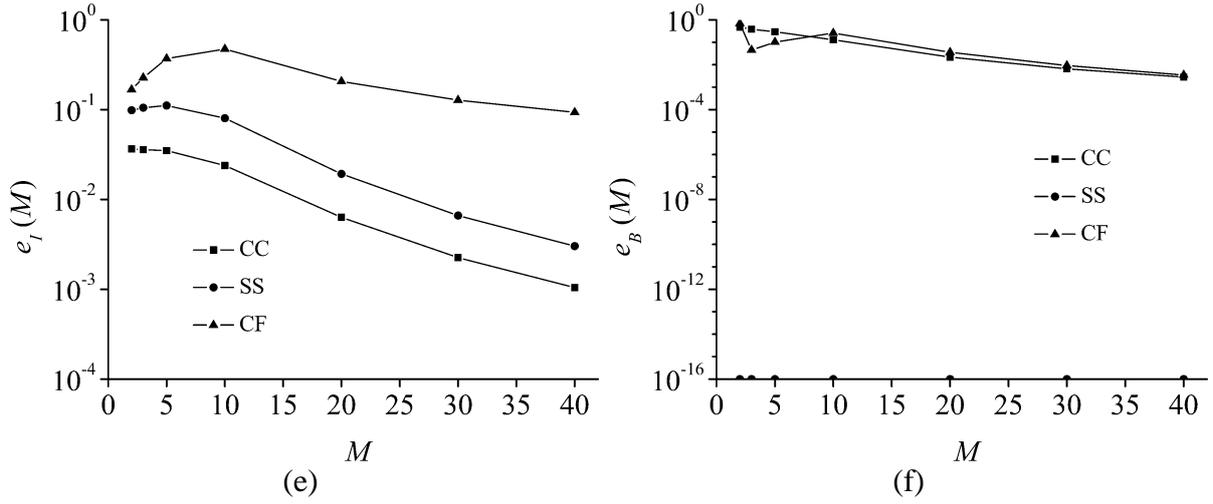

Figure 8: Convergence comparison of the Fourier series multiscale solutions obtained by FCCM, and with different boundary conditions:

(a) $e_I(w)$-$M$ curves, (b) $e_B(w)$-$M$ curves, (c) $e_I(\theta)$-$M$ curves, (d) $e_B(\theta)$-$M$ curves, (e) $e_I(M)$-$M$ curves, (f) $e_B(M)$-$M$ curves.

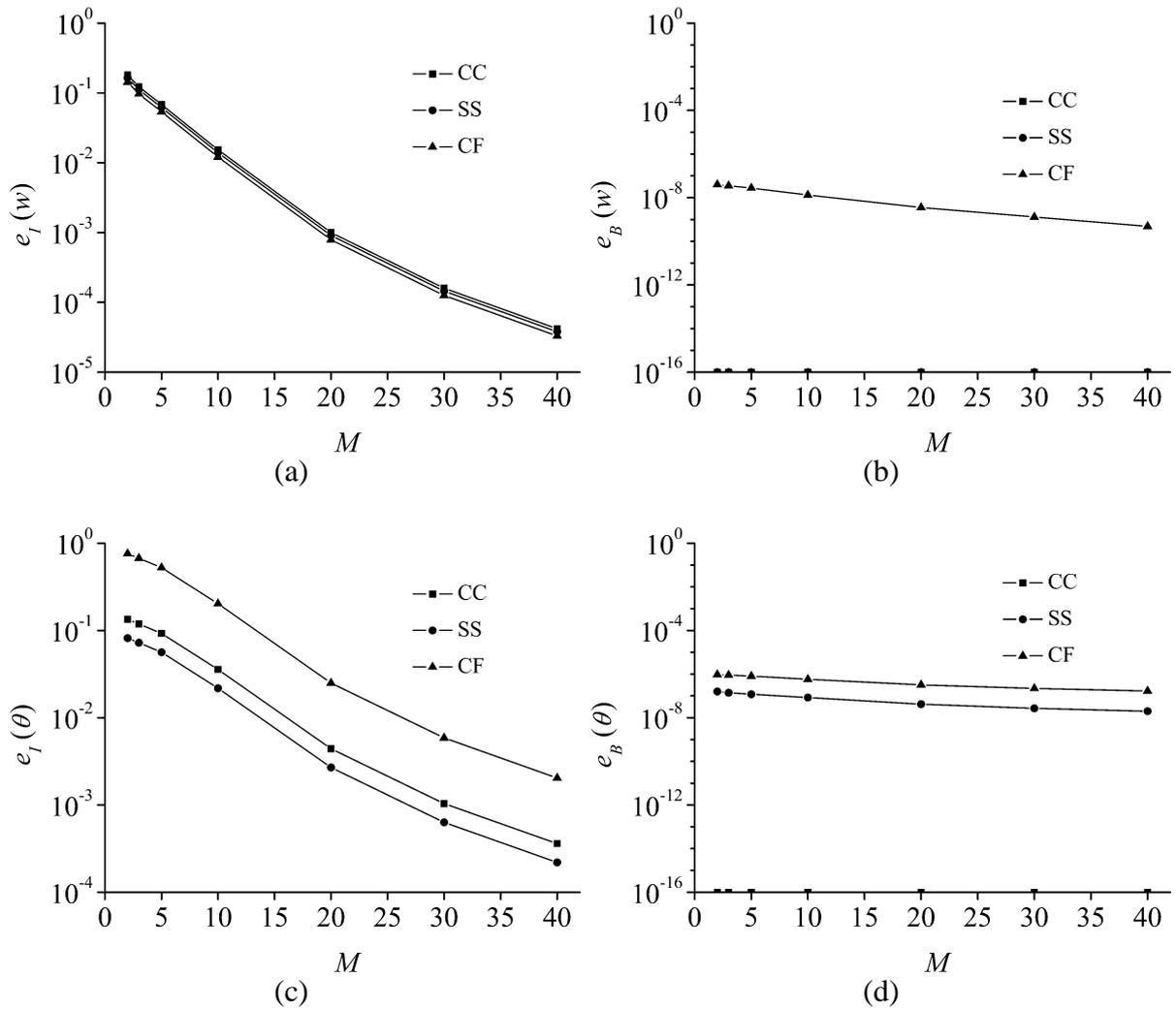

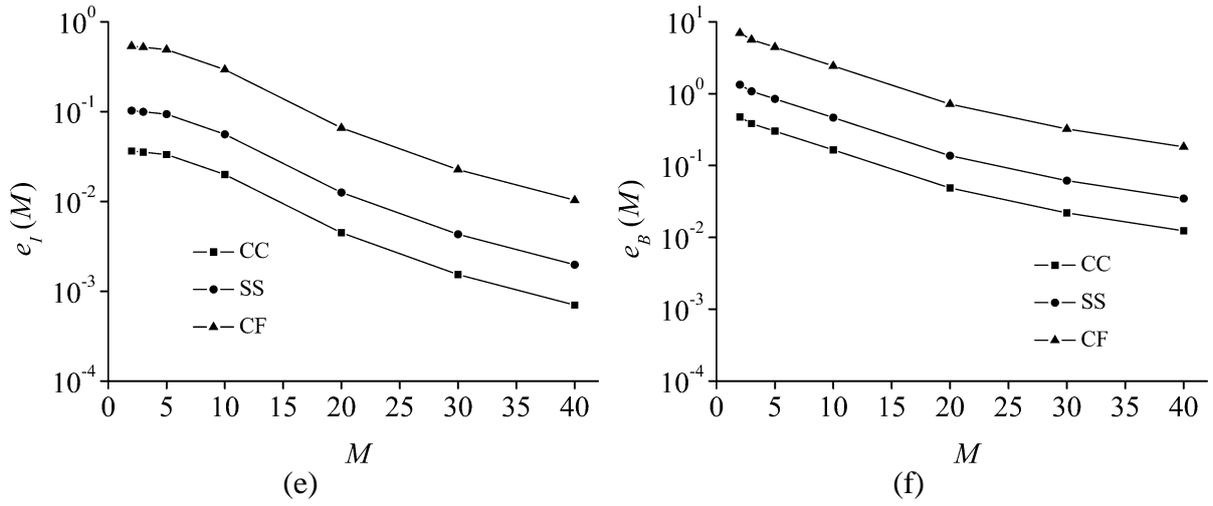

Figure 9: Convergence comparison of the Fourier series multiscale solutions obtained by VM, and with different boundary conditions:
(a) $e_I(w)$-$M$ curves, (b) $e_B(w)$-$M$ curves, (c) $e_I(\theta)$-$M$ curves, (d) $e_B(\theta)$-$M$ curves,
(e) $e_I(M)$-$M$ curves, (f) $e_B(M)$-$M$ curves.

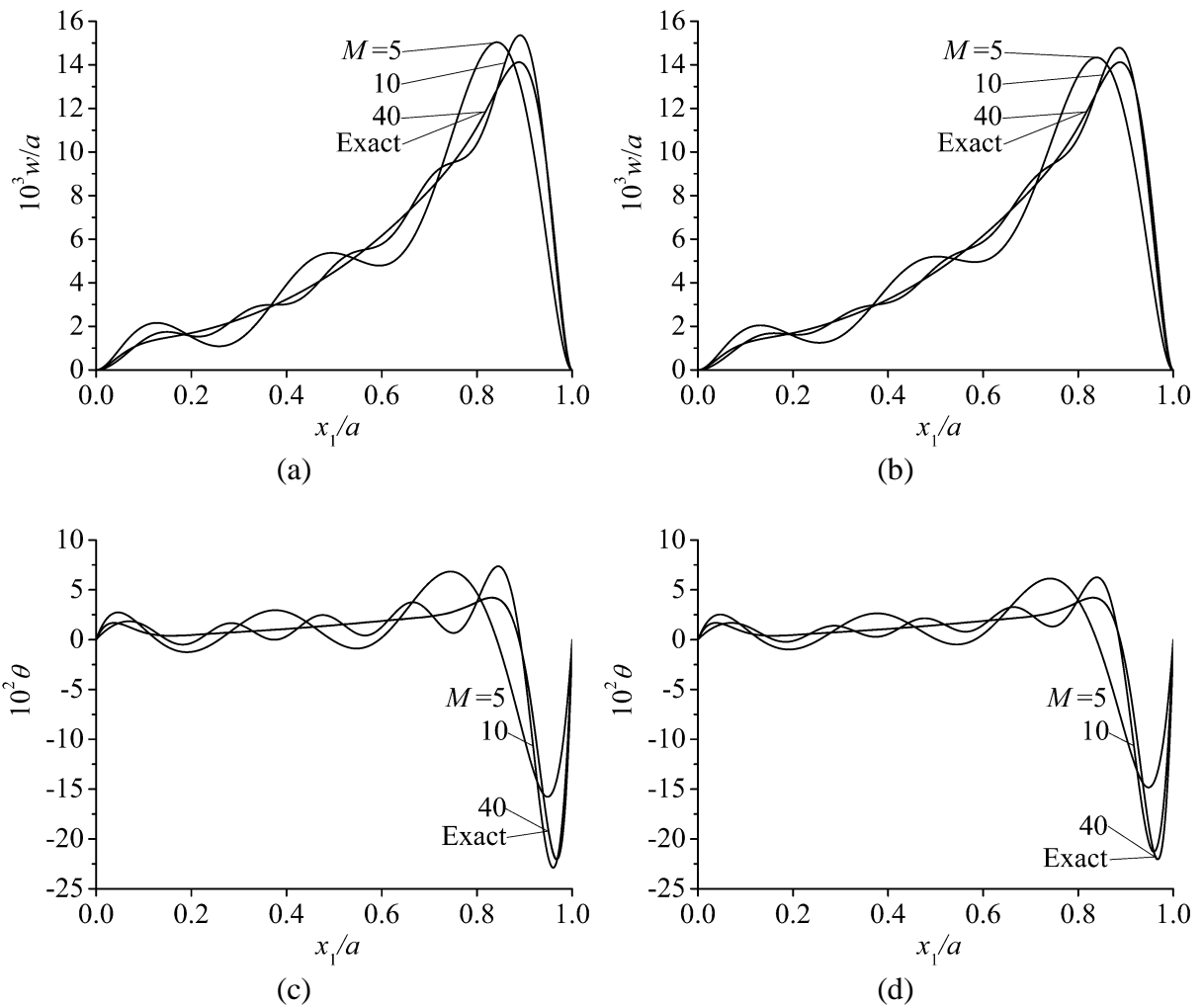

(a)

(b)

(c)

(d)

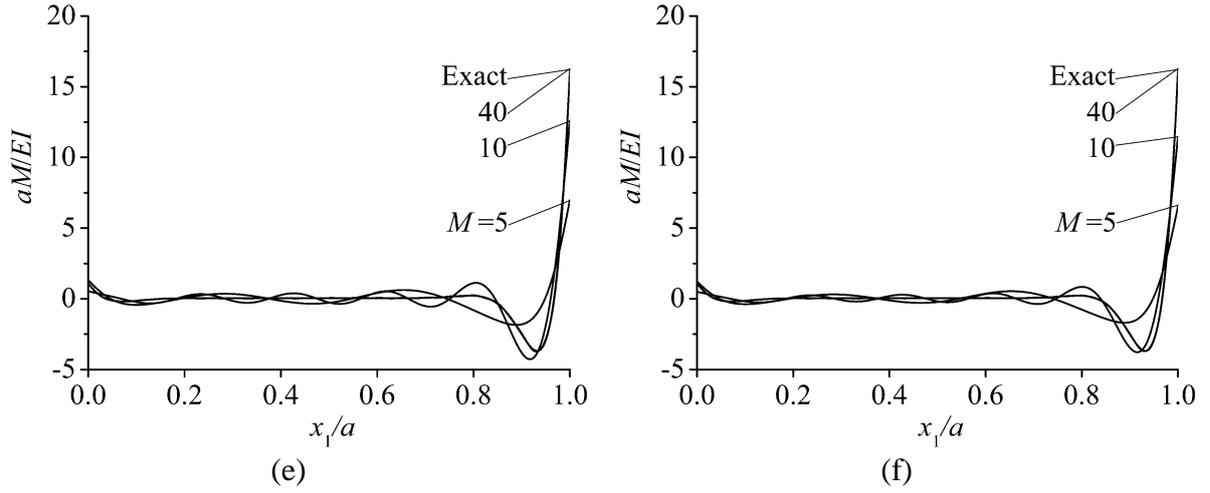

(e)                                                (f)

Figure 10: Convergence comparison of the Fourier series multiscale solutions for the CC boundary condition:
(a) $w(x_1)$-$M$ curves (FCCM), (b) $w(x_1)$-$M$ curves (VM), (c) $\theta(x_1)$-$M$ curves (FCCM),
(d) $\theta(x_1)$-$M$ curves (VM), (e) $M(x_1)$-$M$ curves (FCCM), (f) $M(x_1)$-$M$ curves (VM).

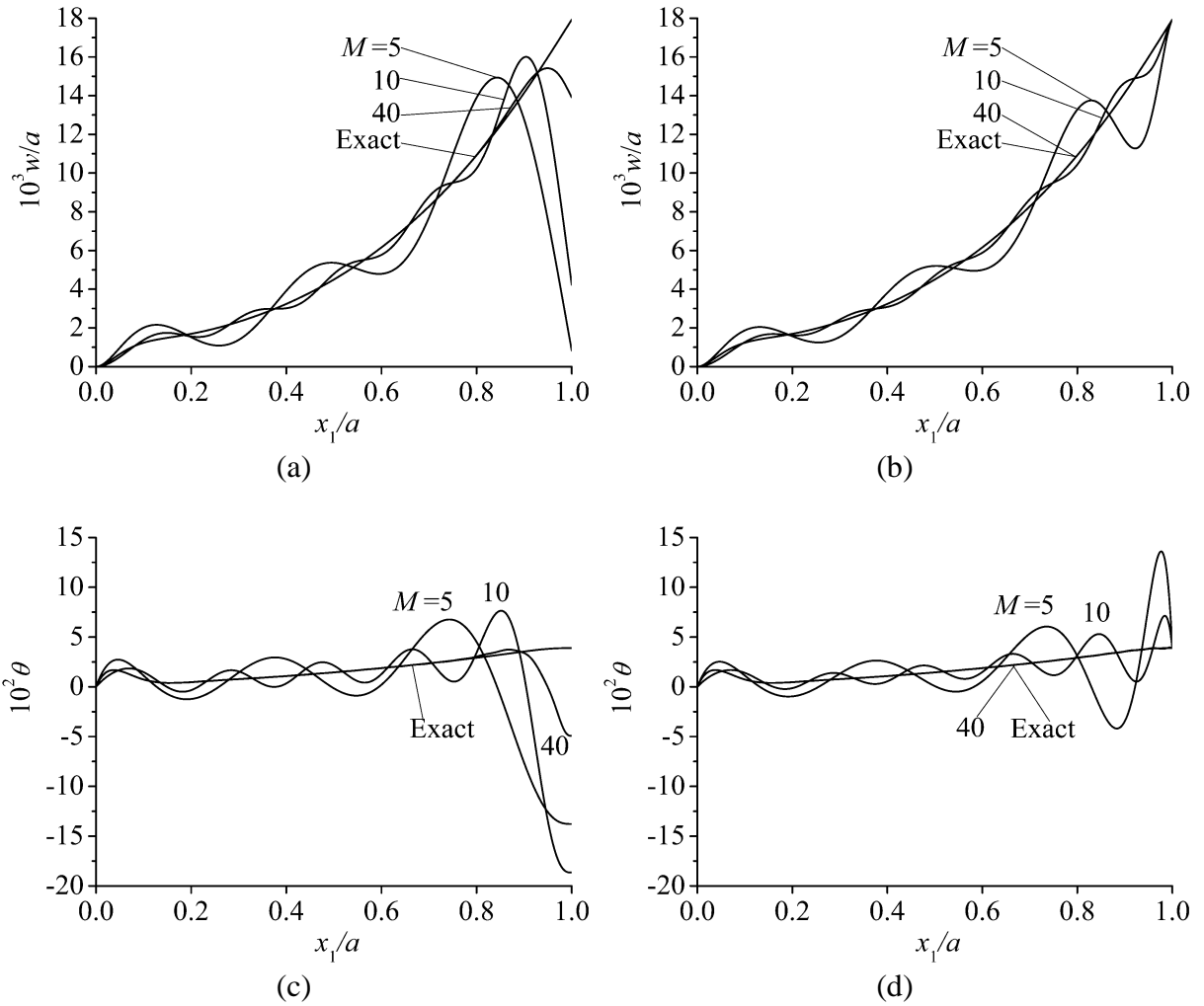

(a)                                                (b)

(c)                                                (d)

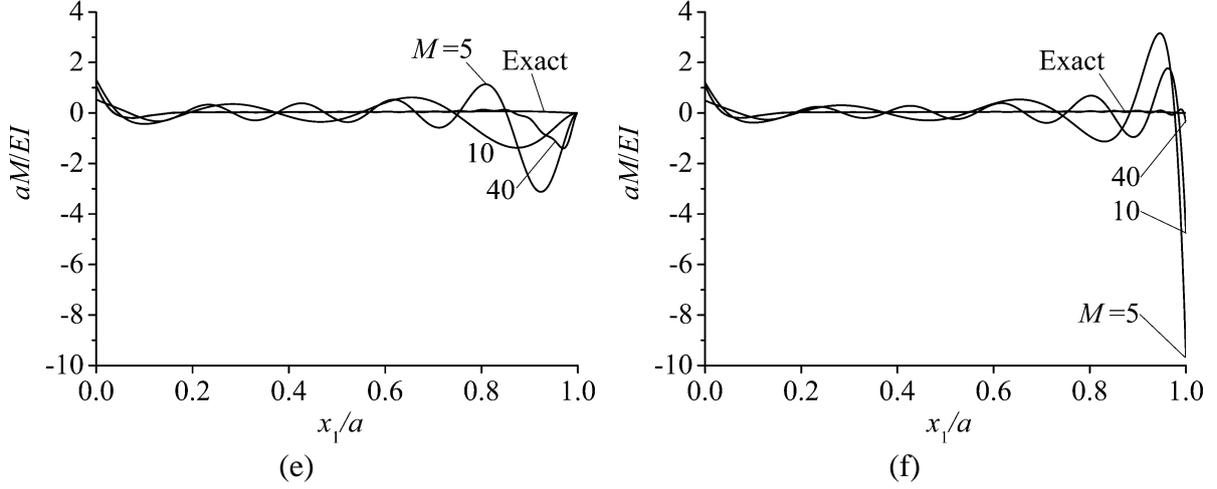

(e)                                             (f)

Figure 11: Convergence comparison of the Fourier series multiscale solutions for the CF boundary condition:
(a) $w(x_1)$-$M$ curves (FCCM), (b) $w(x_1)$-$M$ curves (VM), (c) $\theta(x_1)$-$M$ curves (FCCM),
(d) $\theta(x_1)$-$M$ curves (VM), (e) $M(x_1)$-$M$ curves (FCCM), (f) $M(x_1)$-$M$ curves (VM).

*5.2. Multiscale characteristics*

As shown in Figure 12, the transverse load applied to the beam is of the form
$$q(x_1) = q_0[1000\delta(x_1 - 0.5) + 1000], \quad x_1 \in [0,1], \tag{99}$$
where $q_0 = 1\,\text{N/m}$. We take this quasi-Green's function problem as an example and investigate the multiscale characteristics of elastic beanding of beams on biparametric foundations with varying computational parameters.

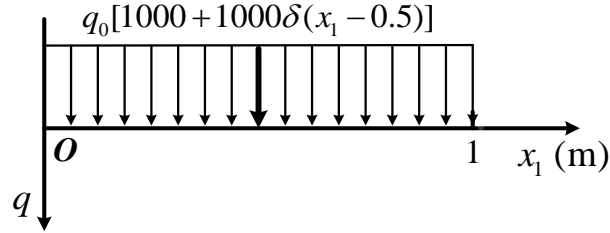

Figure 12: Transverse load applied to the beam.

For the above quasi-Green's function problem, we present in Table 4 the detailed computational scheme.

As shown in Table 4, we pre-specify the computational parameter $k_r$ respectively as 0, $10^2$, $10^4$, and $10^6$ for the elastic bending problem of beam on biparametric foundations and then adjust the computational parameter $G_{pr}$ from 0 to 1000, 2000 and 3000 successively. For these varying computational parameters, we truncate the composite Fourier series with 40 terms, and present in Figure 13 the corresponding transverse relative deflection, $w(x_1)/w(0.5)$, of the beam. It is observed that for smaller values (such as 1.0 and $10^2$) of the parameter $k_r$, the transverse relative deflection undergoes gentle change on the interval $[0,1]$.

As the parameter $k_r$ further increases to $10^4$, a single protrusion appears at action point ($x_1 = 0.5$) of the concentrated load, and the steepness of the protrusion increases significantly with the decrease of the parameter $G_{pr}$. And finally, as the parameter $k_r$ increases to $10^6$, the obvious protrusion appeared at the action point of the concentrated load evolves into a sharp peak. Moreover, for the parameter $G_{pr}=0$, there exist fluctuations on both sides of the peak.

Table 4: Computational scheme for quasi-Green's function problem in elastic beanding of beams on biparametric foundations.

| Description of the problem | | Solution scheme | |
| --- | --- | --- | --- |
| Boundary condition | Computational parameters $(k_r, G_{pr})$ | Order of interpolation algebraical polynomial of load function | Derivation method for discrete equations |
| FF | (1.0, 0) <br> (1.0, 1000) <br> (1.0, 2000) <br> (1.0, 3000) | $N_1^s = 0$ | VM |
| FF | ($10^2$, 0) <br> ($10^2$, 1000) <br> ($10^2$, 2000) <br> ($10^2$, 3000) | $N_1^s = 0$ | VM |
| FF | ($10^4$, 0) <br> ($10^4$, 1000) <br> ($10^4$, 2000) <br> ($10^4$, 3000) | $N_1^s = 0$ | VM |
| FF | ($10^6$, 0) <br> ($10^6$, 1000) <br> ($10^6$, 2000) <br> ($10^6$, 3000) | $N_1^s = 0$ | VM |

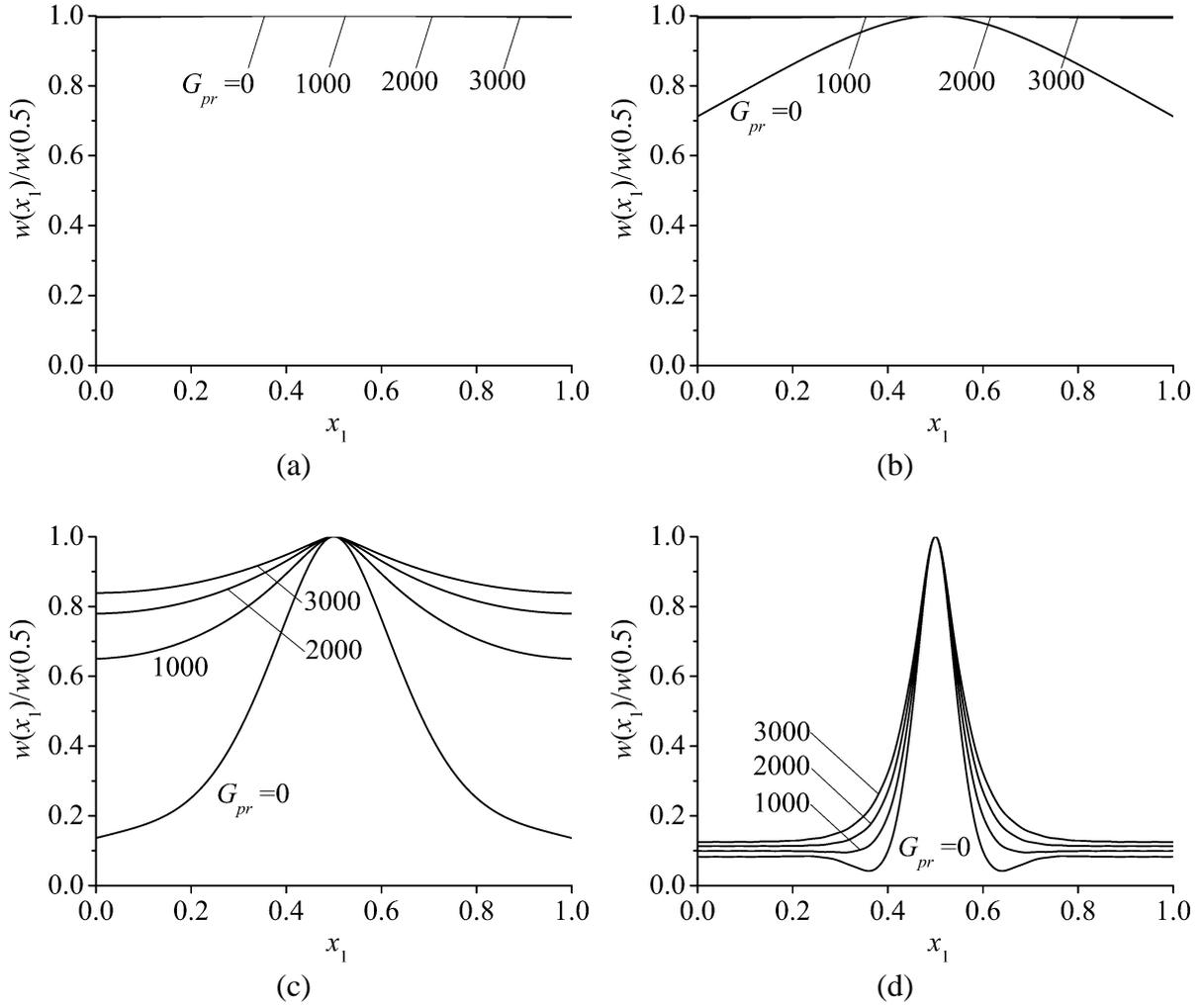

Figure 13: The transverse relative deflection, $w(x_1)/w(0.5)$, of the beam:
(a) $k_r=1.0$, (b) $k_r=10^2$, (c) $k_r=10^4$, (d) $k_r=10^6$.

## 6. Conclusions

    beam on elastic foundations is a common model for several types of engineering structures. In this paper, the usual structural analysis of a Bernoulli-Euler beam on the Pasternak foundation is extended to a multiscale analysis of a fourth order linear differential equation with general boundary conditions and a wide spectrum of model parameters. It is concluded that:

    1. We derive the Fourier series multiscale solution of the bending problem of the Bernoulli-Euler beam on the Pasternak foundation.

    2. We employ the minimum potential energy method, together with the usually used Fourier coefficient comparison method, to derive the final system of equations to be solved for the Fourier series multiscale solution.

    3. We investigate the convergence characteristics of the obtained Fourier series multiscale solution.

    4. We reveal the multiscale characteristics of the bending problem of the Bernoulli-Euler beam on the Pasternak foundation.

The preliminary study on application to bending of Bernoulli-Euler beams on Pasternak foundation verifies the effectiveness of the present Fourier series multiscale method, and provides a new benchmark problem (just as the convection-diffusion-reaction equation) which can be used for persistent improvement in computational performance of other multiscale methods.